\documentclass[a4paper,11pt,twoside]{ps}
\usepackage{amsmath,amssymb,amsthm,graphicx,epic}
\usepackage{rotating}
\usepackage[latin1]{inputenc}
\usepackage[numbers,sort]{natbib}
\let\cite=\citet
\topsep\medskipamount

\newcommand{\asnuc}[1]{\ifcase\csname c@#1\endcsname%
  \or A\or C\or G\or T\else?\fi}
\newcounter{nuc}
\newcounter{dinuc}[nuc]

\newcommand{\carre}[2]{
  \put(-0.5,-0.5){$A(0,0)$}
  \put(-0.5,4.25){$C(0,1)$}
  \put(3.5,-0.5){$T(1,0)$}
  \put(3.5,4.25){$G(1,1)$}
  \multiput(0,0)(#1,0){#2}{\line(0,1){4}}
  \multiput(0,0)(0,#1){#2}{\line(1,0){4}}}

\newcommand{\E}{\mbox{$\mathbb{E}$}}

\newcommand{\PP}{\mbox{$\mathbb{P}$}}
\newcommand{\V}{\mbox{$\mathbb{V}$}}
\renewcommand{\S}[1]{\ensuremath{S\textsc{\lowercase{#1}}}}

\newcommand{\CL}{\mbox{$\mathcal{L}$}}
\newcommand{\CN}{\mbox{$\mathcal{N}$}}
\newcommand{\CA}{\mbox{$\mathcal{A}$}}

\newcommand{\CS}{\mbox{$\mathcal{S}$}}
\newcommand{\CT}{\mbox{$\mathcal{T}$}}
\newcommand{\CX}{\mbox{$\mathcal{X}$}}
\newcommand{\limite}[2]{\mathop{\longrightarrow}
\limits_{\mathrm{#1}}^{\mathrm{#2}}}
\newcommand{\ep}{\mbox{$\varepsilon$}}
\newcommand\1{\leavevmode\hbox{\rm \small1\kern-0.35em\normalsize1}}
\newcommand\ind[1]{\1_{\{#1\}}}
\newcommand\egaldef{\stackrel{\mbox{\upshape\tiny def}}{=}}
\DeclareMathOperator{\dilog}{Li}
\def\DD{\displaystyle}
\begin{document}
\title{Digital Search Trees and Chaos Game Representation}
\author{Peggy Cénac}\address{INRIA Rocquencourt
    and Universit\'{e} Paul Sabatier (Toulouse III) -- INRIA Domaine de
    Voluceau B.P.105 78 153 Le Chesnay Cedex (France)}
\author{Brigitte Chauvin}\address{LAMA, UMR CNRS 8100,
    Bâtiment Fermat, Université de Versailles - Saint-Quentin
    F-78035 Versailles}
\author{Stéphane Ginouillac}\sameaddress{2}
\author{Nicolas Pouyanne}\sameaddress{2}

\begin{abstract}
In this paper, we consider a possible representation of a DNA sequence in a
quaternary tree, in which one can visualize repetitions of subwords
(seen as suffixes of subsequences).
The CGR-tree turns a sequence of letters into a Digital Search Tree (DST),
obtained from the suffixes of the reversed sequence.
Several results are known concerning the height, the insertion depth for DST
built from independent successive random sequences having the same
distribution.
Here the successive inserted words are strongly dependent.
We give the asymptotic behaviour of the insertion depth and
the length of branches for the CGR-tree obtained from the suffixes of a
reversed i.i.d.\@ or Markovian sequence.
This behaviour turns out to be at first order the same one as in the case of
independent words.
As a by-product, asymptotic results on the length of longest runs in a
Markovian sequence are obtained.
\end{abstract}
\begin{resume}

La repr\'esentation d\'efinie ici est une repr\'esentation possible de
s\'equence d'ADN dans un arbre quaternaire dont la construction permet de
visualiser les r\'ep\'etitions de suffixes.
\`A partir d'une s\'equence de lettres, on construit un arbre digital de
recherche (\emph{Digital Search Tree}) sur l'ensemble des suffixes de la
s\'equence invers\'ee.
Des r\'esultats sur la hauteur et la profondeur d'insertion ont \'et\'e
\'etablis lorsque les s\'equences \`a placer dans l'arbre sont
ind\'ependantes les unes des autres.
Ici les mots \`a ins\'erer sont fortement d\'ependants.
On donne le comportement asymptotique de la profondeur d'insertion et
de la longueur des branches pour un arbre obtenu \`a partir des suffixes d'une
s\'equence i.i.d.\@ ou markovienne retourn\'ee.
Au premier ordre, cette asymptotique est la m\^eme que dans le cas o\`u les
mots ins\'er\'es sont ind\'ependants.
De plus, certains r\'esultats peuvent aussi s'interpr\'eter comme des
r\'esultats de convergence sur les longueurs de plus longues r\'ep\'etitions
d'une lettre dans une s\'equence Markovienne.
\end{resume}


\maketitle  
\section{Introduction}\label{sec:intro}
In the last years, DNA has been represented by means of several methods
in order to make pattern visualization easier and to detect local or global
similarities (see for instance \cite{Roy}).
The \emph{Chaos Game Representation} (CGR) provides both a graphical
representation and a storage tool. From a sequence in a finite alphabet,
CGR defines a trajectory in a bounded subset of $\xR ^d$ that keeps all
statistical properties of the sequence.
\cite{Jeffrey} was the first to apply this iterative method to DNA sequences.
\cite{CFL, PeggytestCGR} study the CGR with
an extension of word-counting based methods of analysis.
In this context, sequences are made of $4$ nucleotides named A (adenine),
C (cytosine), G (guanine) and T (thymine).

\begin{figure}
\label{Fi:compte}
\begin{center}
\begin{center}
\begin{picture}(11,5)(-2.5,-2.5)
\put(-2,-2){\carre{2}{3}}
\put(-1.2,-1.2){\S{A}}
\put(-1.2,0.8){\S{C}}
\put(0.8,-1.2){\S{T}}
\put(0.8,0.8){\S{G}}
\put(4,-2){\carre{1}{5}}
\put(4.2,-1.6){\small{\S{AA}}}
\put(4.2,-0.6){\small{\S{CA}}}
\put(4.2,0.4){\small{\S{AC}}}
\put(4.2,1.4){\small{\S{CC}}}
\put(5.2,-1.6){\small{\S{TA}}}
\put(5.2,-0.6){\small{\S{GA}}}
\put(5.2,0.4){\small{\S{TC}}}
\put(5.2,1.4){\small{\S{GC}}}
\put(6.2,-1.6){\small{\S{AT}}}
\put(6.2,-0.6){\small{\S{CT}}}
\put(6.2,0.4){\small{\S{AG}}}
\put(6.2,1.4){\small{\S{CG}}}
\put(7.2,-1.6){\small{\S{TT}}}
\put(7.2,-0.6){\small{\S{GT}}}
\put(7.2,0.4){\small{\S{TG}}}
\put(7.2,1.4){\small{\S{GG}}}
\end{picture}
\end{center}
\medskip
\includegraphics[width=7cm,height=7cm]{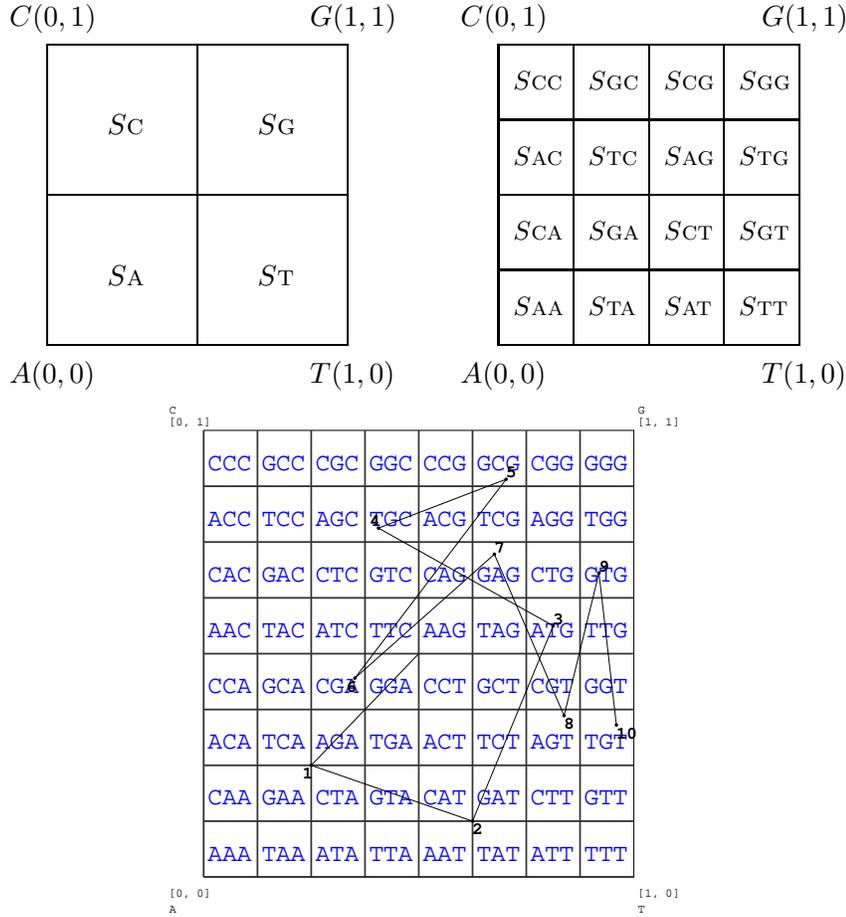}
\end{center}
\caption{Chaos Game Representation of the first $10$ nucleotides of the
  \emph{E. Coli} threonine gene thrA\@: ATGCGAGTGT\@. The coordinates
  for each nucleotide are calculated recursively using $(0.5, 0.5)$ as
  starting position. The sequence is read from left to right. Point
  number $3$ corresponds to the first $3$-letter word $ATG$. It is
  located in the corresponding quadrant. The second $3$-letter word
  $TGC$ corresponds to point $4$ and so on.}
\end{figure}

The CGR of a sequence $U_{1} \ldots U_{n}\ldots$ of letters $U_{n}$ from a
finite alphabet $\CA$ is the sequence $(\CX _{n})_{n\geq 0}$ of points in an
appropriate compact subset $S$ of $\xR ^d$ defined by
\[\left\{
\begin{array}{l}
 \CX _{0} \in S \\
 \CX _{n+1}=\theta\bigl(\CX _{n}+\ell_{U_{n+1}}\bigr),
\end{array}
\right.
 \]
where $\theta$ is a real parameter ($0<\theta<1$), each letter $u\in\CA$ being
assigned to a given point $\ell_u\in S$.
In the particular case of Jeffrey's representation, $\CA =\{ A,C,G,T\}$ is the
set of nucleotides,
$S=[0,1]^2$ is the unit square.
Each letter is placed at a vertex as follows:
\[
  \ell_A=(0,0),\quad \ell_C=(0,1), \quad \ell_G=(1,1),\quad\ell_T=(1,0),
\]
$\theta=\frac{1}{2}$ and the first point $\CX _{0}$ is the center of the
square.
Then, iteratively, the point $\CX _{n+1}$ is the middle of the segment between
$\CX _{n}$ and the square's vertex $\ell_{U_{n+1}}$:
\[\CX _{n+1}=\frac{\CX _{n}+\ell_{U_{n+1}}}{2},\]
or, equivalently,
\[
\CX _{n}=\sum_{k=1}^{n}\frac{\ell_{U_{k}}}{2^{n-k+1}}+\frac{\CX _{0}}{2^{n}}.
\]
Figure \ref{Fi:compte} 
represents the construction of the word
ATGCGAGTGT.

With each deterministic word $w=u_{1}\ldots u_{n}$, we associate the
half-opened subsquare $Sw$ defined by the formula
\[
Sw\egaldef \sum_{k=1}^{n}\frac{\ell_{u_{k}}}{2^{n-k+1}}
+\frac{1}{2^n}[0,1[^{2};
\]
it has center $\sum_{k=1}^{n}\ell_{u_{k}}/2^{n-k+1}+\CX _{0}/2^{n}$
and side $1/2^n$.
For a given random or deterministic sequence $U_{1} \ldots U_{n}\ldots$,
for any word $w$ and any $n\geq |w|$ (the notation $|w|$ stands for the number
of letters in $w$),
counting the number of points $(\CX _{i})_{1\leq i\leq n}$ that
belong to the subsquare $Sw$ is tantamount to counting the number
of occurences of $w$ as a subword of $U_{1} \ldots U_{n}$.
Indeed, all successive words from the sequence having
$w$ as a suffix are represented in $Sw$.
See Figure~\ref{Fi:compte} for an example with three-letter subwords.
This provides tables of word frequencies (see \cite{Goldman}).
One can generalize it to any subdivision of the unit square;
when the number of subsquares is not a power of $4$, the table of word
frequencies defines a counting of words with noninteger length
(see \cite{Almeida}). 

The following property of the CGR is important:
\it
the value of any $\CX _n$
contains the historical information of the whole sequence
$\CX _{1},\ldots \CX _{n}$.
\rm
Indeed, notice first that, by construction, $\CX _n\in Su$ with $U_n=u$;
the whole sequence is now given by the inductive formula
$\CX _{n-1}=2\CX _n-\ell_{U_n}$.

We define a representation of a random DNA sequence $U=(U_n)_{n\geq 1}$
as a random
quaternary tree, the {\it CGR-tree}, in which one can visualize repetitions
of subwords.
We adopt the classical order $(A,C,G,T)$ on letters.
Let $\CT$ be the complete infinite $4$-ary tree;
each node of $\CT$ has four branches corresponding to letters $(A,C,G,T)$
that are ordered in the same way.
The CGR-tree of $U$ is an increasing sequence
$\CT_{1} \subset \CT_{2} \ldots \subset \CT_{n} \subset\ldots$
of finite subtrees of $\CT$, each $\CT_{n}$ having $n$ nodes.
The $\CT_{n}$'s are built by successively inserting the {\it reversed prefixes}
\begin{equation}
\label{reversed}
W(n)=U_n\ldots U_1
\end{equation}
as follows in the complete infinite tree.
First letter $W(1)=U_{1}$ is inserted in the complete infinite tree at level
$1$, {\it i.e.} just under the root, at the node that corresponds to the letter
$U_1$.
Inductively, the insertion of the word $W(n)=U_n\ldots U_1$ is made as follows:
try to insert it at level $1$ at the node $\CN$ that corresponds to
the letter $U_{n}$.
If this node $\CN$ is vacant, insert $W(n)$ at $\CN$ ;
if $\CN$ is not vacant, try to insert $W(n)$ in the subtree having $\CN$ as
a root, at the node that corresponds to the letter $U_{n-1}$, and so on.
One repeats this operation until the node at level $k$ that corresponds
to letter $U_{n-k+1}$ is vacant;
word $W(n)$ is then inserted at that node.

We complete our construction by labelling the $n$-th inserted node with the
word $W(n)$.
One readily obtains this way the process of a digital search tree (DST),
as stated in the following proposition.

Figure~2
shows the very first steps of construction of the tree
that corresponds to any sequence that begins with $GAGCACAGTGGAAGGG$.
The insertion of this complete 16-letter prefix is represented in
Figure~3.
In these figures, each node has been labelled by its order of insertion to make
the example more readable.

\begin{figure}[p]
\label{Fi:pasApas}
\def\hei{9.5truecm}
\def\wid{8.3truecm}
\begin{center}
\includegraphics[width=\wid ,height=\hei ]{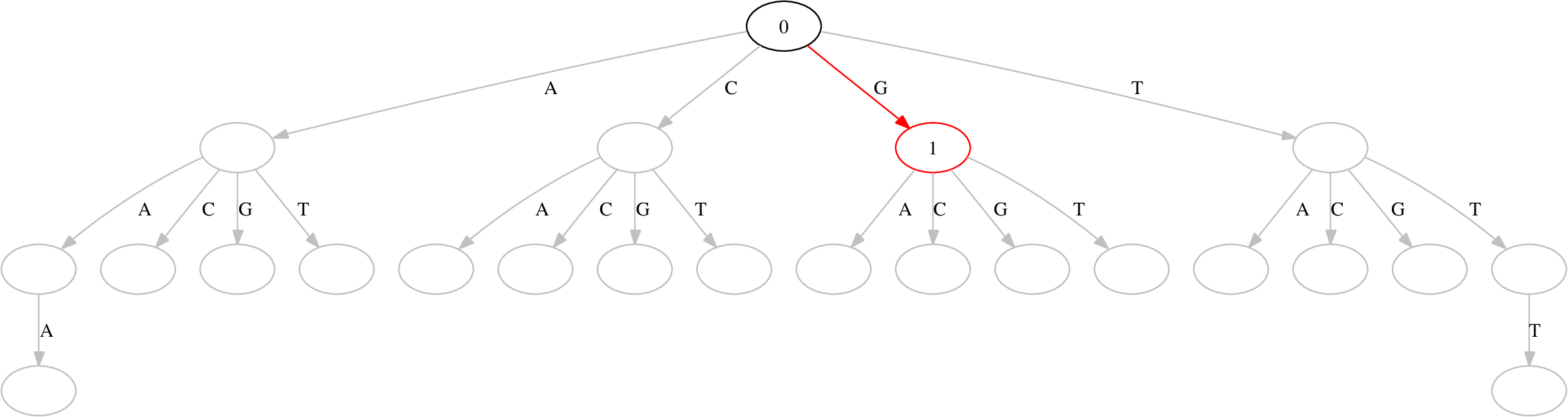}
\hskip 2truemm
\includegraphics[width=\wid ,height=\hei ]{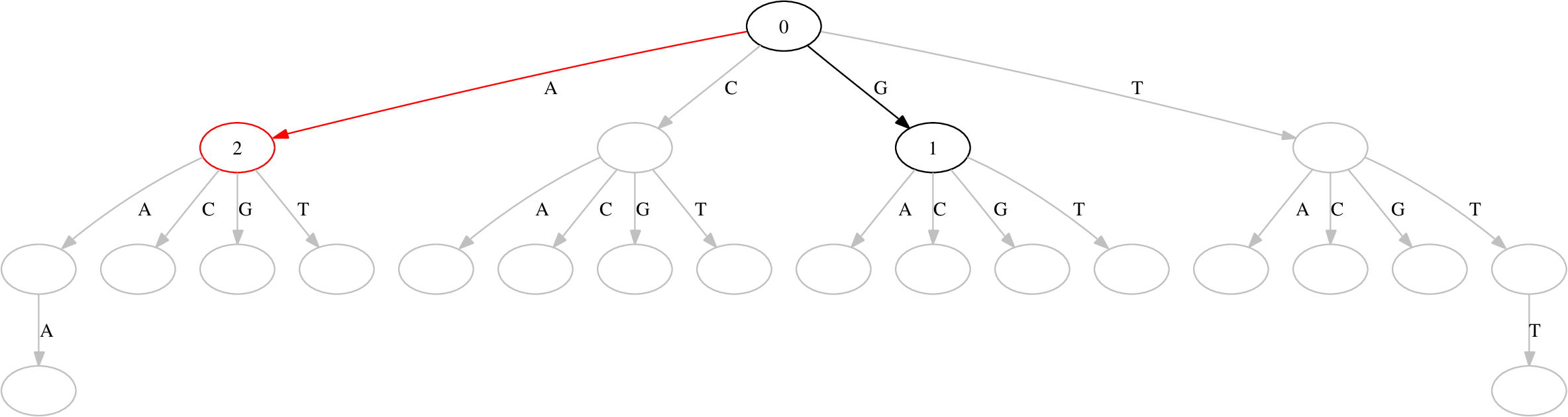}
\end{center}
\begin{center}
\includegraphics[width=\wid ,height=\hei ]{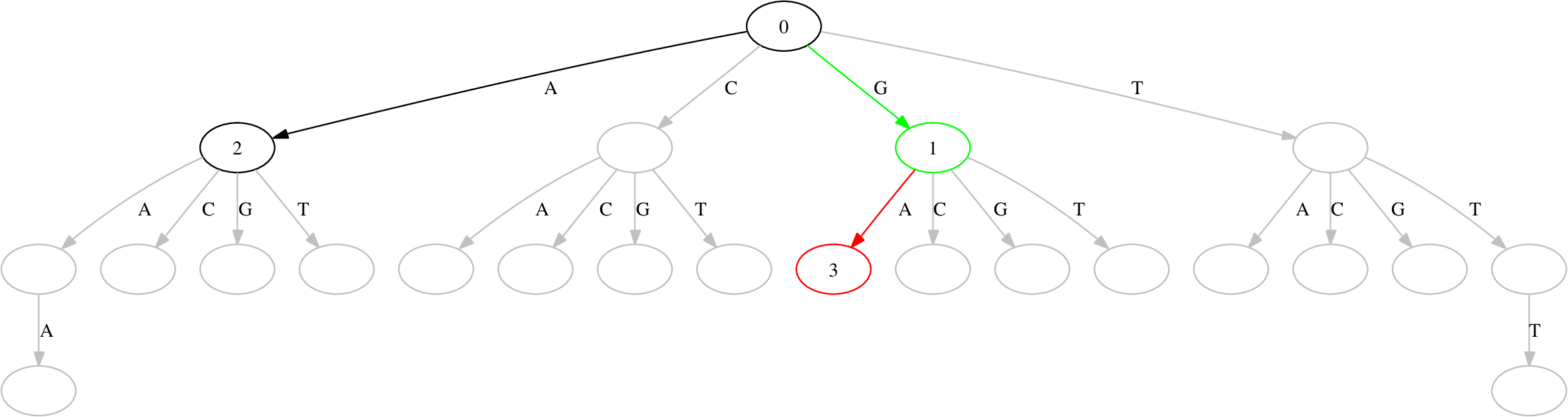}
\hskip 2truemm
\includegraphics[width=\wid ,height=\hei ]{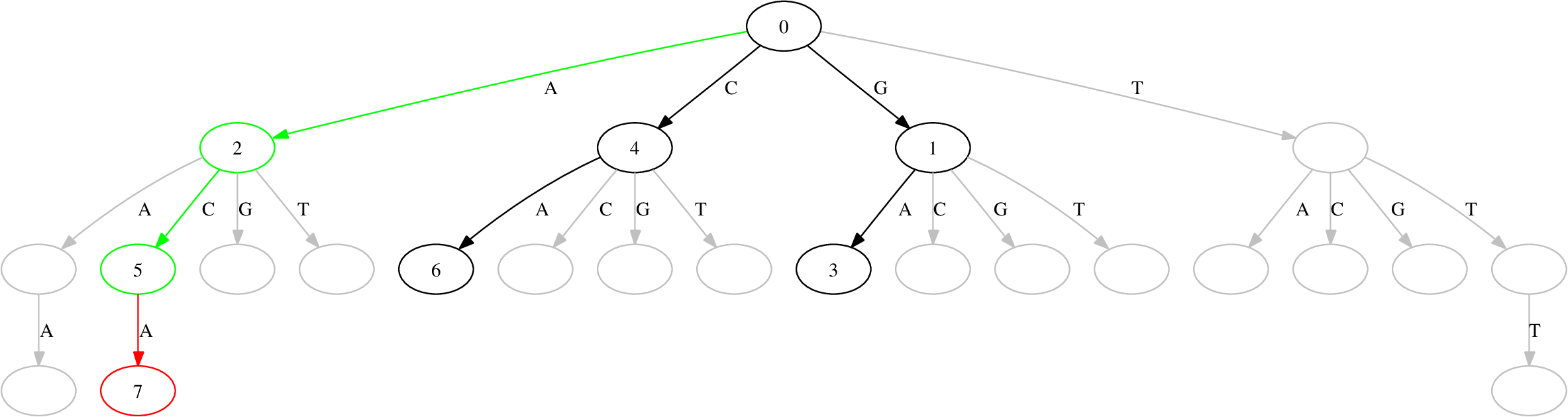}
\end{center}
\caption{Insertion of a sequence $GAGCACAGTGGAAGGG...$ in its CGR-tree:
first, second, third and seventh steps.}
\end{figure}

\begin{prpstn}
The CGR-tree of a random sequence $U=U_1U_2\dots$ is a digital search tree, 
obtained by insertion in a quaternary tree of the successive reversed prefixes
$U_1$, $U_2U_1$, $U_3U_2U_1$, \dots of the sequence.
\end{prpstn}


\begin{figure}[p]
\label{Fi:arbre}
\begin{center}
\includegraphics[width=7cm,height=6cm]{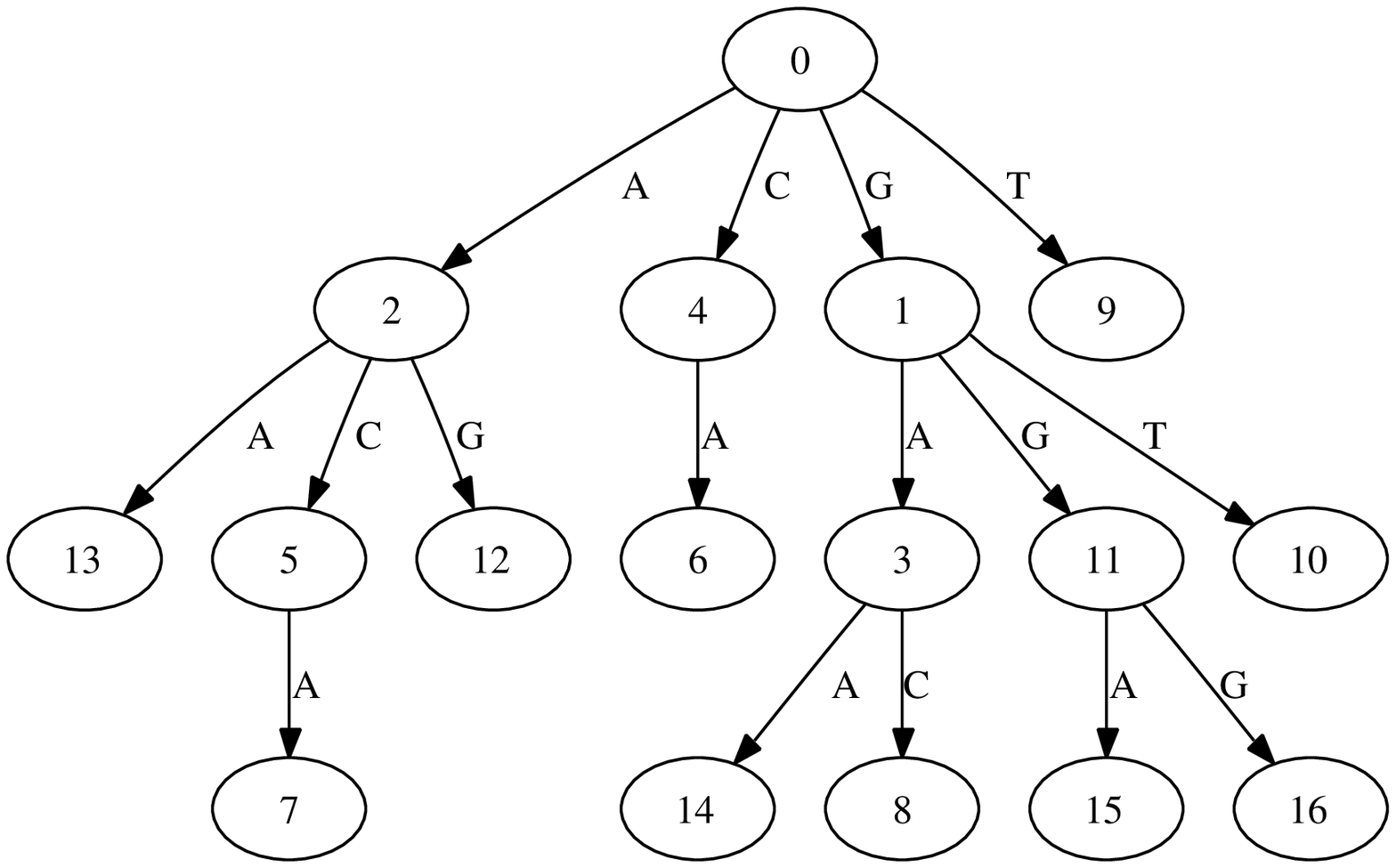}
\quad
\includegraphics[width=5cm,height=5cm]{dessinnorm.ps}
\end{center}
\caption{Representation of $16$ nucleotides of
\emph{Mus Musculus}\@ GAGCACAGTGGAAGGG\@ in the CGR-tree (on the left)
and in the ``historyless representation'' (on the right).}
\end{figure}

\medskip
The main results of our paper are the following convergence results,
the random sequence $U$ being supposed to be Markovian.
If $\ell _n$ and $\CL _n$ denote respectively the length of the shortest
and of the longest branch of the CGR-tree, then {\it $\ell _n/\ln n$ and
$\CL _n/\ln n$ converge almost surely to some constants}
(Theorem \ref{cvps}).
Moreover, if $D_n$ denotes the insertion depth and if $M_n$ is the length
of a uniformly chosen random path, then {\it $D_n/\ln n$ and $M_n/\ln n$
converge in probability to a common constant}
(Theorem \ref{thm}).

\medskip
\begin{rmrk}
A given CGR-tree without its labels (i.e.\@ a given shape of tree) is
equivalent to a list of words in the sequence without their order.
More precisely, one can associate with a shape of
CGR-tree, a representation in the unit square as described below.
With any node of the tree (which is in bijection with a word
$w=W_{1}\ldots W_{d}$), we associate the center of the corresponding square
$Sw$,
\[
\CX _{w}\egaldef
\sum_{k=1}^{d}\frac{\ell_{W_{k}}}{2^{d-k+1}}+\frac{\CX _{0}}{2^{d}}.
\] 
For example, Figure~3
shows this ``\emph{historyless representation}'' for the
word $GAGCACAGTGGAAGGG$. Moreover
Figure~4
enables us to qualitatively compare the original and
the historyless representations on an example.
\end{rmrk}
\begin{figure}[p]
\label{Fi:compare}
\begin{center}
\includegraphics[width=0.45\textwidth]{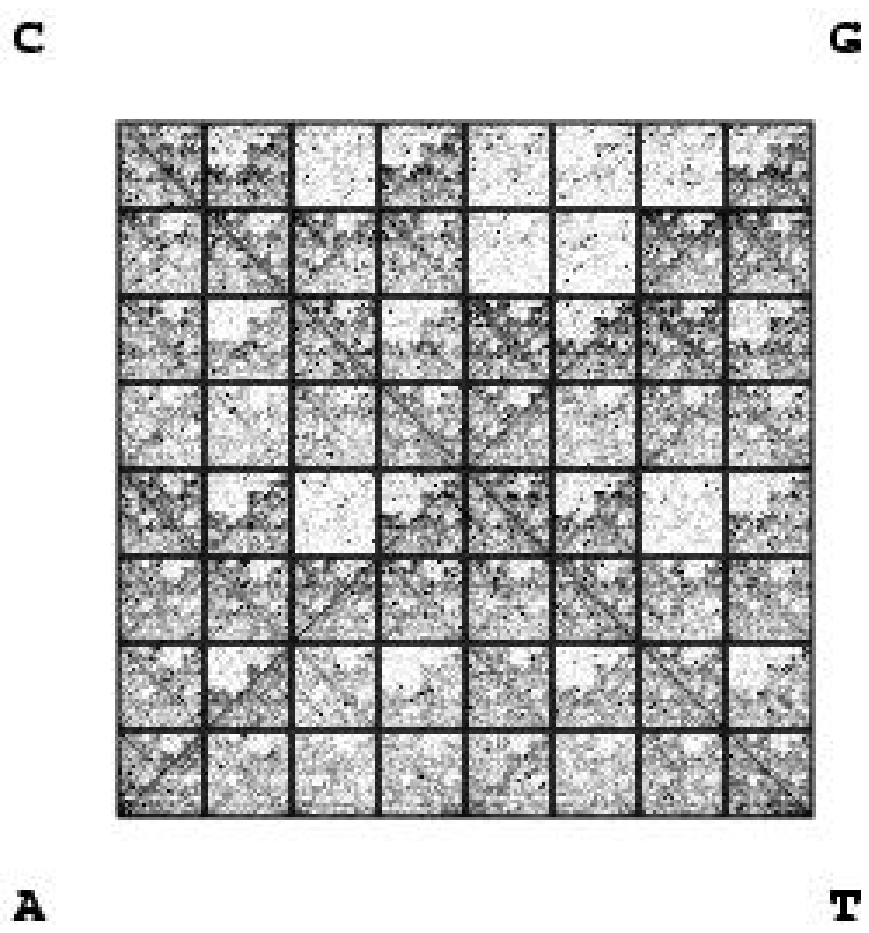}
\hfill
\includegraphics[width=0.45\textwidth]{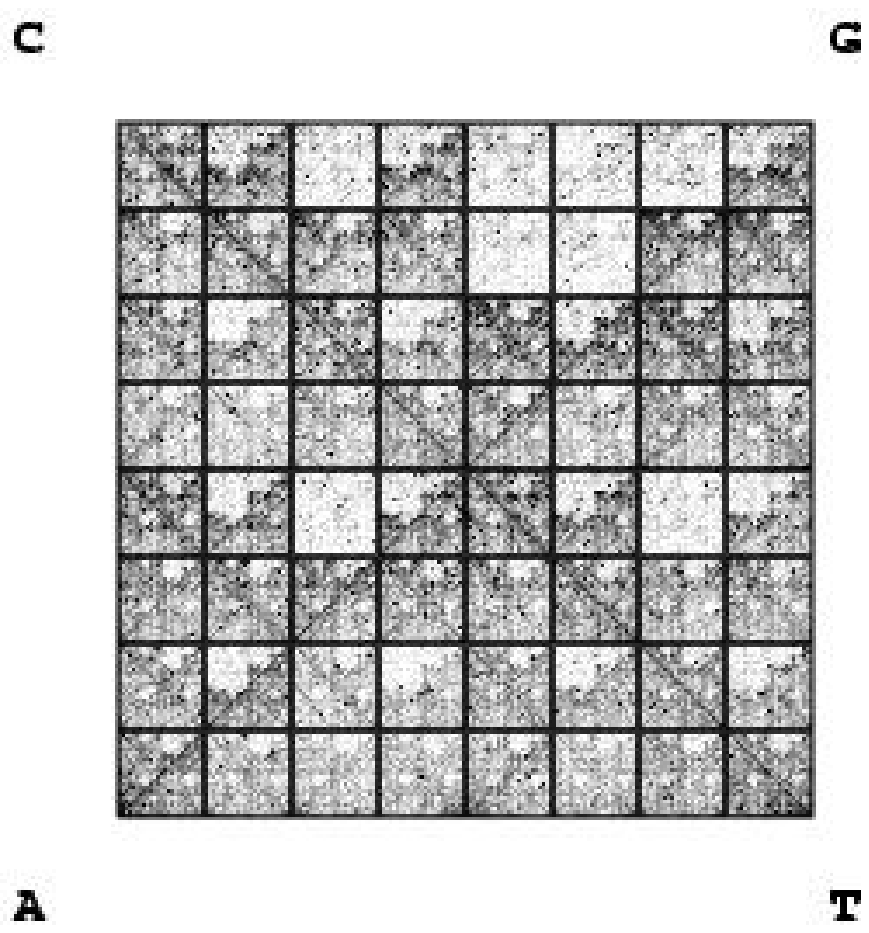}
\caption{Chaos Game Representation (on the left) and historyless
representation (on the right) of the first $400 000$ nucleotides of
Chromosome 2 of \emph{Homo Sapiens}.}
\end{center}
\end{figure}

Several results are known (see chap. 6 in \cite{Mahmoud}), concerning the
height, the insertion depth and the profile for DST obtained from
{\it independent} successive sequences, having the same distribution. 
It is far from our situation where the successive inserted words are strongly
dependent from each other.
Various results concerning the so-called Bernoulli model
(binary trees, independent sequences
and the two letters have the same probability 1/2 of appearance)
can be found in \cite{Mahmoud}.
\cite{Aldous} prove by embedding in continuous time, that the height satisfies
$H_{n}-\log_{2}n \rightarrow 0$ in probability.
Also \cite{Drmota} proves that the height of such DSTs is concentrated:
$\E[H_{n}-\E(H_{n})]^{L}$ is asymptotically bounded for any $L>0$.

For DST constructed from independent sequences on an $m$-letter alphabet
with nonsymmetric (i.e. non equal probabilities on the letters) i.i.d
or Markovian sources, \cite{Pittel} gets several results on the insertion
depth and on the height.
Despite the independence of the sequences, Pittel's work seems to be the
closest to ours, and some parts of our proofs are inspired by it.

Some proofs in the sequel use classical results on the distribution of word
occurences in a random sequence of letters
(independent or Markovian sequences).
\cite{Blom} give the generating function of the first occurence of a word
for i.i.d.\@ sequences, based on a recurrence relation on the probabilities.
This result is extended to Markovian sequences by \cite{Daudin-Robin}. 
Several studies in this domain are based on generating functions, for example
\cite{Regnier}, \cite{Schbath}, \cite{SP}.
Nonetheless, other approaches are considered:
one of the more general techniques is the Markov chain embedding method
introduced by \cite{Fu} and further developped by \cite{Fu-Koutras, Koutras}.
A martingale approach (see \cite{Li, Gerber-Li, Williams}) is an alternative to
the Markov chain embedding method to solve problems around \cite{Penney} Game.
These two approaches are compared in \cite{Pozdnyakov}.
Whatever method one uses, the distribution of the first occurence of a word
strongly depends on its overlapping structure.
This dependence is at the core of our proofs.

As a by-product, our results yield asymptotic properties on the length of the
longest run, which is a natural object of study.
In i.i.d.\@ and symmetric sequences, \cite{Revesz} establish almost sure
results about the growth of the longest run.
These results are extended to Markov chains in \cite{Samarova}, and
\cite{Gordon} show that the probabilistic behaviour of the length of the
longest run is closely approximated by that of the maximum of some i.i.d.\@
exponential random variables.

The paper is organized as follows.
In Section~\ref{sec:def} we establish the assumptions and notations we use
throughout.
Section \ref{sec:longueur} is devoted to almost sure convergence of the
shortest and the longest branches in CGR-trees.
In Section~\ref{sec:profondeur} asymptotic behaviour of the insertion depth is
studied.
An appendix deals separately with the domain of definition of the generating
function of a certain waiting time related to the overlapping structure of
words.


\section{Assumptions and notations}\label{sec:def}

In all the sequel, the sequence $U=U_1\dots U_n \dots $ is supposed to be a
Markov chain of order $1$, with transition matrix  $Q$ and invariant measure
$p$ as initial distribution.

For any deterministic infinite sequence $s$, let us denote by $s^{(n)}$ the
word formed by the $n$ first letters of $s$, that is to say
$s^{(n)}\egaldef s_{1}\dots s_{n}$, where $s_{i}$ is the $i$-th
letter of $s$. 
The measure $p$ is extended to reversed words the following way:
$p(s^{(n)})\egaldef \PP(U_{1}=s_{n}, \ldots, U_{n}=s_{1})$.
The need for reversing the word $s^{(n)}$ comes from the construction of the
CGR-tree which is based on reversed sequences (\ref{reversed}).

We define the constants
\begin{eqnarray*}
h_+&\egaldef&\lim_{n\to +\infty}\frac{1}{n}\max\Bigl\{
\ln\biggl(\frac{1}{p\bigl(s^{(n)}\bigr)}\biggr) ,~ p\bigl(s^{(n)}
\bigr)>0\Bigr\},\\
h_-&\egaldef&\lim_{n\to +\infty}\frac{1}{n}\min\Bigl\{
\ln\biggl(\frac{1}{p\bigl(s^{(n)}\bigr)}\biggr) ,~p\bigl(s^{(n)}
\bigr)>0\Bigr\},\\
h&\egaldef&\lim_{n\to +\infty}\frac{1}{n}
\E\Bigl[\ln\Bigl(\frac{1}{p\bigl(U^{(n)}\bigr)}\Bigr)\Bigr].
\end{eqnarray*}
Due to an argument of sub-additivity (see \cite{Pittel}), these limits
are well defined (in fact, in a  more general than Markovian
sequences framework).
Moreover, Pittel proves the existence of two infinite sequences denoted here
by $s_{+}$ and  $s_{-}$ such that 
\begin{equation}
\label{h+-}
h_+=\lim_{n \to \infty}\frac{1}{n}
\ln\biggl(\frac{1}{p\bigl(s_{+}^{(n)}\bigr)}\biggr),
\quad \mbox{and}\quad
h_-=\lim_{n \to \infty}\frac{1}{n}
\ln\biggl(\frac{1}{p\bigl(s_{-}^{(n)}\bigr)}\biggr).
\end{equation}
For any $n \geq 1$, the notation $\CT_{n}\egaldef \CT_{n}(W)$ stands for
the finite tree with
$n$ nodes (without counting the root), built from the first $n$ sequences
$W(1), \ldots, W(n)$, which are the successive reversed prefixes of the
sequence $(U_n)_n$, as defined by~(\ref{reversed}).
$\CT_{0}$ denotes the tree reduced to the root.
In particular, the random trees are increasing:
$ \CT_{0}\subset \CT_{1}\ldots \subset \CT_{n} \subset \ldots\subset \CT$. 

Let us define $\ell_{n}$ (resp. $\CL_{n}$) as the length of the shortest
(resp. the longest) path from the root to a feasible external node of the tree
$\CT_{n}(w)$.
Moreover,  $D_{n}$ denotes the insertion depth of $W(n)$ in $\CT_{n-1}$ to
build $\CT_{n}$.
Finally $M_{n}$ is the length of a path of $\CT_{n}$, randomly and
uniformly chosen in the $n$ possible paths. 

The following random variables play a key role in the proofs. For the
sake of precision, let us recall that $s$ is deterministic, the randomness is
uniquely due to the generation of the sequence $U$.
First we define for any infinite sequence $s$ and for any $n \geq 0$,
\begin{equation}
\label{defXn}
X_{n}(s)\egaldef \left \{
\begin{array}{l}
0 \ \mbox{if}\ s_{1}\ \mbox{is not in}\ \CT_{n}\\
\max \{k ~\mbox{such~that}\ s ^{(k)}\ \mbox{is already inserted in}
\ \CT_{n}\}. \\
\end{array}
\right. 
\end{equation}
Notice that $X_0(s) =0$.
Every infinite sequence corresponds to a branch of the infinite tree $\CT$
(root at level $0$, node that corresponds to $s_1$ at level $1$,
node that corresponds to $s_2$ at level $2$, {\it etc.});
the random variable $X_{n}(s)$ is the length of the branch associated with $s$
in the tree $\CT_{n}$.
For any $k \geq 0$,
$T_{k}(s)$ denotes the size of the first tree where $s^{(k)}$ is inserted:
\[
T_{k}(s)\egaldef \min\{n,~X_{n}(s)=k\}
\]
(notice that $T_0(s) =0$).

These two variables are in duality in the following sense:
one has equality of the events
\begin{equation}
\label{duality}
\{ X_{n}(s) \geq k \}=\{ T_{k}(s) \leq n \}
\end{equation}
and consequently, $\{ T_{k}(s) = n \} \subset \{ X_{n}(s) = k \}$ since
$X_n(s)-X_{n-1}(s)\in\{ 0,1\}$.

In our example of Figures~2
and~3,
the drawn random
sequence is $GAGCACAGTGGAAGGG\dots$
If one takes a deterministic sequence $s$ such that $s^{(3)}=ACA$, then
$X_0(s)=X_1(s)=0$, $X_2(s)=X_3(s)=X_4(s)=1$, $X_5(s)=X_6(s)=2$ and
$X_k(s)=3$ for $7\leq k\leq 18$.
The first three values of $T_k(s)$ are consequently
$T_1(s)=2$, $T_2(s)=5$, $T_3(s)=7$.

Moreover, the random variable $T_{k}(s)$ can be decomposed as follows,
\begin{equation}
\label{decompcle}
T_{k}(s)=\sum_{r=1}^{k}Z_{r}(s),
\end{equation}
where $Z_{r}(s)\egaldef T_{r}(s)-T_{r-1}(s)$ is the number of letters to read
before the branch that corresponds to $s$ increases by $1$.
In what follows,
$Z_{r}(s)$ can be viewed as the waiting time $n$ of the first occurence of
$s^{(r)}$ in the sequence 
\[
\ldots U_{n+T_{r-1}(s)}U_{n-1+T_{r-1}(s)}\ldots U_{1+T_{r-1}(s)}s^{(r-1)},
\]
i.e. $Z_{r}(s)$ can also be defined as
\[
Z_{r}(s)=\min\{n \geq 1,~ U_{n+T_{r-1}(s)}\ldots
U_{n+T_{r-1}(s)-r+1}=s_{1} \ldots s_{r} \}.
\]
Because of the Markovianity of the model, the random variables $Z_{r}(s)$
are independent.

Let us then introduce $Y_r(s)$ as being the waiting time of the first
occurence of $s^{(r)}$ in the sequence
\[\ldots U_{n+T_{r-1}(s)}U_{n-1+T_{r-1}(s)}\ldots U_{1+T_{r-1}(s)},\]
that is to say
\[
Y_{r}(s)=\min\{n \geq r,~U_{n+T_{r-1}(s)}\ldots
U_{n+T_{r-1}(s)-r+1}=s_{1} \ldots s_{r} \}.
\]
One has readily the inequality $Z_{r}(s) \leq Y_r(s)$.
More precisely, if the word $s^{(r)}$ is inserted in the sequence before
time $T_{r-1}(s)+r$, there is some overlapping between prefixes of $s^{(r-1)}$
and suffixes of $s^{(r)}$.
See Figure~5
for an example where $r=6$ and
$s_1s_2s_3=s_4s_5s_6$.
Actually, variables $Z_{r}(s)$ and $Y_r(s)$ are related by
\[
Z_{r}(s)=\ind{Z_{r}(s)<r}Z_{r}(s)+\ind{Z_{r}(s)\geq r}Y_{r}(s).
\]
Since the sequence $(U_{n})_{n\geq 1}$ is stationary, the conditional
distribution of $Y_r(s)$ given $T_{r-1}(s)$ is the distribution of the first
occurence of the word $s^{(r)}$ in the realization of a Markov chain of order
$1$, whose transition matrix is $Q$ and whose initial distribution is its
invariant measure.
In particular the conditional distribution of $Y_r(s)$ given $T_{r-1}(s)$ is
independent of $T_{r-1}(s)$.
\setlength{\unitlength}{0.8cm}
\begin{figure}
\label{Fi:overlap}
\begin{center}
\begin{picture}(17,5)(0,0)
\put(0,0){\line(10,0){12}}
\put(0,0.6){\line(10,0){12}}

\put(-1,0.8){\line(12,0){18}}
\put(-1,1.45){\line(12,0){18}}

\put(0,0){\line(0,1){0.6}}
\put(2,0){\line(0,1){0.6}}
\put(4,0){\line(0,1){0.6}}
\put(6,0){\line(0,1){0.6}}
\put(8,0){\line(0,1){0.6}}
\put(10,0){\line(0,1){0.6}}
\put(12,0){\line(0,1){0.6}}

\put(0,0.8){\line(0,1){0.65}}
\put(2,0.8){\line(0,1){0.65}}
\put(4,0.8){\line(0,1){0.65}}
\put(6,0.8){\line(0,1){0.65}}
\put(8,0.8){\line(0,1){0.65}}
\put(10,0.8){\line(0,1){0.65}}
\put(12,0.8){\line(0,1){0.65}}
\put(14,0.8){\line(0,1){0.65}}
\put(16,0.8){\line(0,1){0.65}}

\put(0.9,0.2){$s_1$}
\put(2.9,0.2){$s_2$}
\put(4.9,0.2){$s_3$}
\put(6.9,0.2){$s_4$}
\put(8.9,0.2){$s_5$}
\put(10.9,0.2){$s_6$}
\put(0.2,1){$U_{3+T_{5}(s)}$}
\put(2.2,1){$U_{2+T_{5}(s)}$}
\put(4.2,1){$U_{1+T_{5}(s)}$}
\put(6.9,1){\bf$s_1$}
\put(8.9,1){\bf$s_2$}
\put(10.9,1){\bf$s_3$}
\put(12.9,1){\bf$s_4$}
\put(14.9,1){\bf$s_5$}
\end{picture}
\end{center}
\caption{How overlapping intervenes in $Z_{r}(s)$' definition.
In this example, one takes $r=6$.
In the random sequence, prefix $s^{(6)}$ can occur starting from
$U_{3+T_{5}(s)}$ only if $s_{1}s_{2}s_{3}=s_{4}s_{5}s_{6}$.}
\end{figure}
\setlength{\unitlength}{1cm}

The generating function $\Phi(s^{(r)},t)\egaldef \E[t^{Y_r(s)}]$ is given by
\cite{Daudin-Robin}:
\begin{equation}
\label{gen}
\Phi (s^{(r)},t)=\Bigl(\gamma_{r}(t)+(1-t)\delta_{r}(t^{-1})\Bigr)^{-1},
\end{equation}
where the functions $\gamma$ and $\delta$ are respectively defined as
\begin{equation}
\label{gammadelta}
\gamma_{r}(t)\egaldef \frac{1-t}{tp\bigl(s_{r}\bigr)}
\sum_{m \geq 1}Q^{m}(s_{1}, s_{r})t^{m},
\quad
\delta_{r}(t^{-1})\egaldef \sum_{m=1}^r
\frac
{\ind{s_{r} \ldots s_{r-m+1}=s_{m} \ldots s_{1}}}
{t^{m}p\bigl(s^{(m)}\bigr)},
\end{equation}
and where $Q^{m}(u,v)$ denotes the transition probability from $u$ to $v$
in $m$ steps.

\begin{rmrk}
In the particular case when the sequence of nucleotides $(U_n)_{n\geq 1}$ is
supposed to be independent and identically distributed according to the non
degenerated law $(p_{A}, p_{C}, p_{G}, p_{T})$, the transition probability
$Q^{m}(s_{1},s_{r})$ is equal to $p(s_{r})$, and hence $\gamma_{r}(t)=1$.
\end{rmrk}

\begin{prpstn}
\label{generatrice}
\begin{itemize}
\item[(i)]
The generating function of $Y_{r}(s)$ defined by~(\ref{gen})
has a ray of convergence $\geq 1+\kappa p\bigl(s^{(r)}\bigr)$
where $\kappa$ is a positive constant independent of $r$ and $s$.
\item[(ii)]
Let $\gamma$ denote the second largest eigenvalue of the transition matrix
$Q$.
For all $t \in ]-\gamma^{-1},\gamma^{-1}[$, 
\begin{equation}
\label{inegGamma}
\bigl|\gamma_{r}(t)-1\bigr| \leq \frac{|1-t|}{1-\gamma |t|} \kappa',
\end{equation}
where $\kappa'$ is some positive constant independent of $r$ and $s$
(if $\gamma =0$ or if the sequence is i.i.d.\@, we adopt the convention
$\gamma^{-1}=+\infty$ so that the result remains valid).
\end{itemize}
\end{prpstn}

\begin{proof}
The proof of Proposition~\ref{generatrice} is given in Appendix~\ref{sec:A}.
\end{proof}

\section{Length of the branches}\label{sec:longueur}


In this section we are concerned with the asymptotic behaviour of the length
$\ell_{n}$ (resp. $\CL_{n}$) of the shortest (resp. longest) branch of the
CGR-tree.

\begin{thrm}
\label{cvps}
\[\frac{\ell_{n}}{\ln n}\limite{n \to \infty}{a.s.} \frac{1}{h_+},\quad
\mbox{and} \quad 
\frac{\CL_{n}}{\ln n}\limite{n \to \infty}{a.s.}\frac{1}{h_-}.\]
\end{thrm}

According to the definition of $X_{n}(s)$, the lengths $\ell_{n}$ and
$\CL_{n}$ are  functions of $X_{n}$:
\begin{equation}
\label{longueurs}
\ell_{n}=\min_{s \in \CA ^{\xN}}X_{n}(s), \quad \mbox{and} \quad
\CL_{n}=\max_{s \in \CA ^{\xN}}X_{n}(s).
\end{equation}
The following key lemma gives an asymptotic result on $X_{n}(s)$, under
suitable assumptions on $s$.
Our proof of Theorem \ref{cvps} is based on it.

\begin{lmm}
\label{lemme}
Let $s$ be such that there exists 
\begin{equation}
\label{assumption}
\lim_{n\longrightarrow +\infty}
\frac{1}{n}\ln\biggl(\frac{1}{p\bigl(s ^{(n)}\bigr)}\biggr)\egaldef h(s)>0.
\end{equation}
Then
\[
\frac{X_{n}(s)}{\ln n} \limite{n \to \infty}{a.s.}\frac{1}{h(s)}.
\]
\end{lmm}

\begin{rmrk}
Let $\tilde v \egaldef vv\ldots$ consist of repetitions of a letter $v$. 
Then $X_{n}(\tilde v)$ is the length of the branch associated with $\tilde v$
in $\CT_{n}$.
For such a sequence (and exclusively for them) the random variable
$Y_{k}(\tilde v)$ is equal to  $T_{k}(\tilde v)$.
Consequently $X_{n}(\tilde v)$ is the length of the longest run of '$v$'
in $U_{1} \ldots U_{n}$.
When $(U_n)_{n\geq 1}$ is a sequence of i.i.d.\@ trials,
\cite{Petrov, Renyi, Revesz} showed that
\[
\frac{X_{n}(\tilde v)}{\ln n}
\limite{n \to \infty}{a.s.}\frac{1}{\ln\frac{1}{p}},
\]
where $p\egaldef \PP(U_{i}=v)$.
This convergence result is a particular case of Lemma~\ref{lemme}.  
\end{rmrk}

\begin{figure}[p]
\label{Fi:simulations1}
\begin{center}
\includegraphics[height=8.5truecm,width=15truecm]{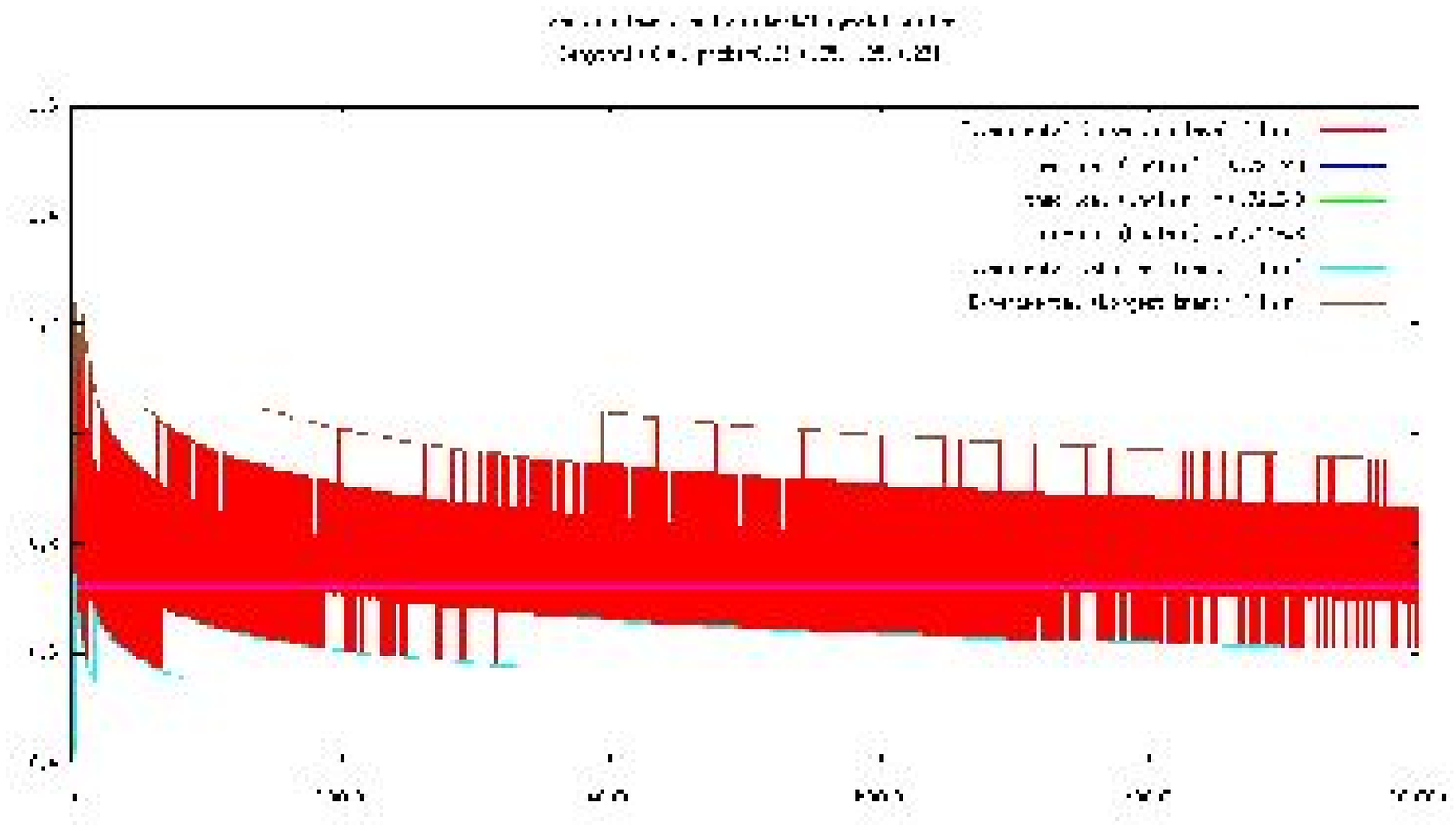}
\includegraphics[height=8.5truecm,width=15truecm]{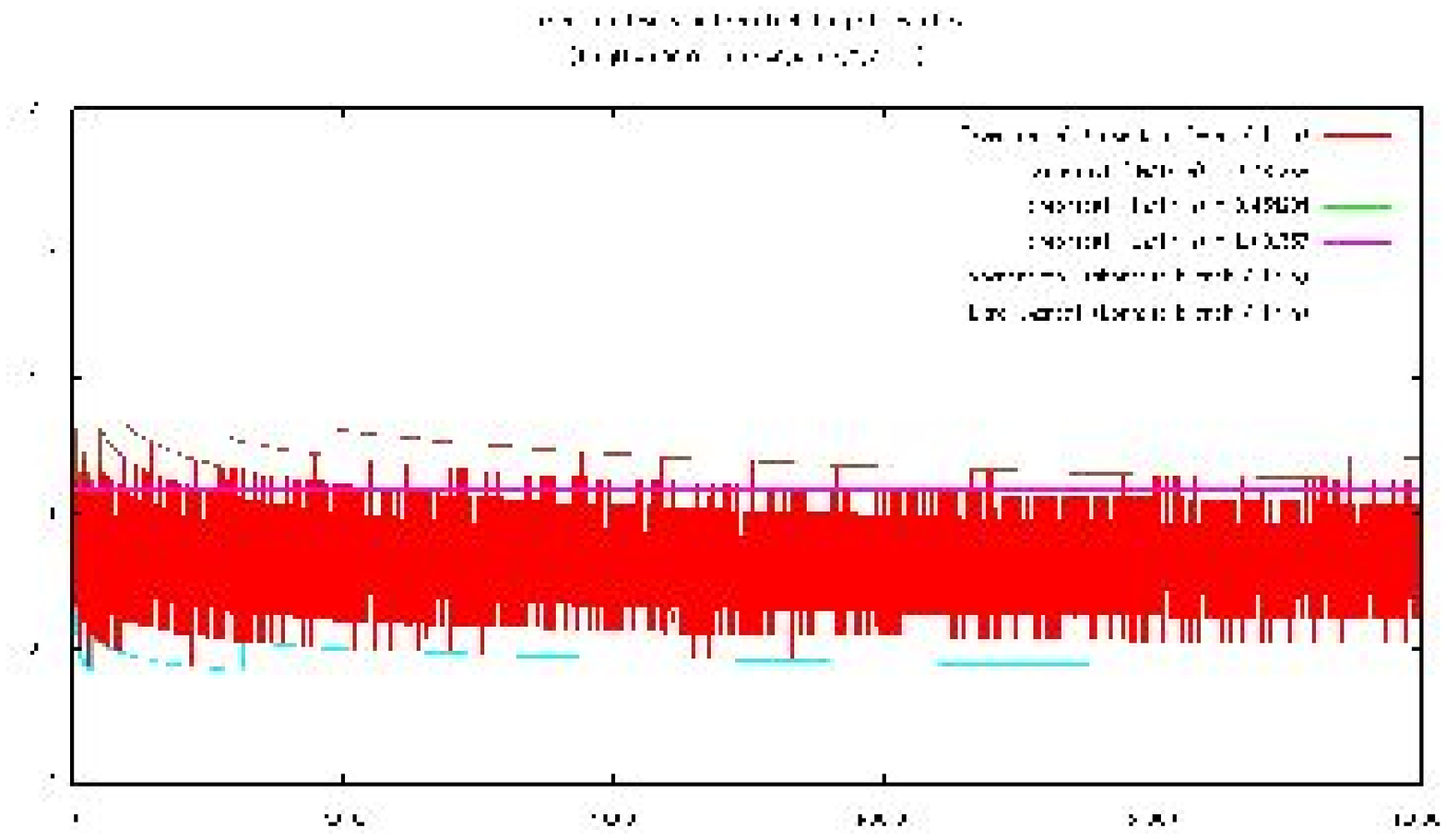}
\caption{Simulations of two random sequences.
On the first graphic, letters of the sequence are i.i.d.~and equally likely
distributed; on the second one, i.i.d.~letters have probabilities
$(p_A,p_C,p_G,p_T)=(0.4,0.3,0.2,0.1)$.
On the $x$-axis, number $n$ of inserted letters;
on the $y$-axis, normalized insertion depth $D_n/\ln n$ (oscillating curve),
lengths of the shortest and of the longest branch (regular ``under'' and
``upper envelops'').
The horizontal lines correspond to the constant limits of these three random
variables (on the first graph, these three limits have the same value).}
\end{center}
\end{figure}
\begin{figure}[p]
\label{Fi:simulations2}
\begin{center}
\includegraphics[width=15truecm,height=15truecm]{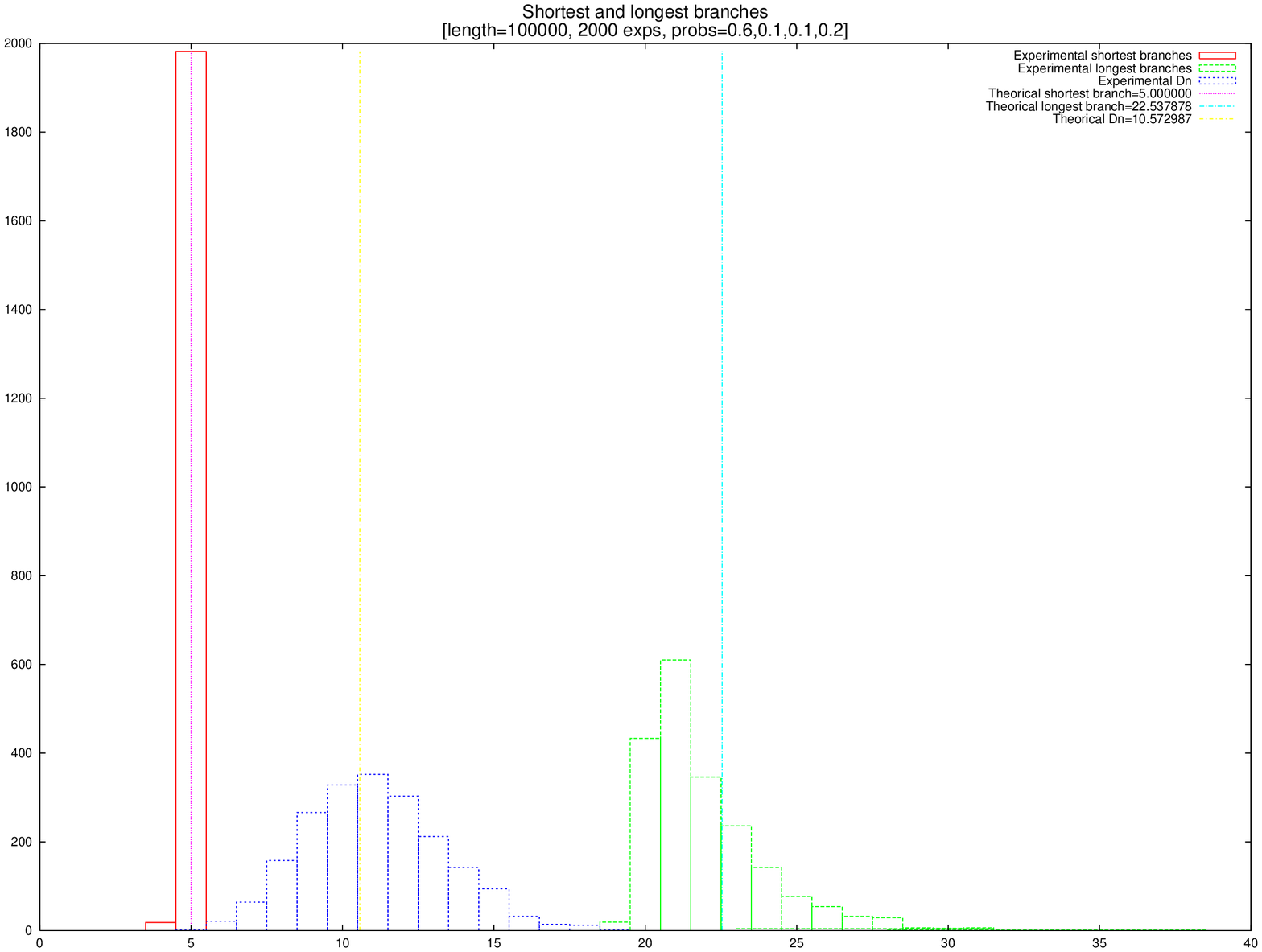}
\caption{Simulations of $2000$ sequences of $100,000$ i.i.d.~letters.
On the left, histogram of shortest branches;
in the middle, histogram of insertion depth of the last inserted word;
on the right, histogram of longest branches.
Vertical lines are their expected values, namely $\ln (10^5)\times \ell$ where
$\ell$ respectively equals the limit of $\ell_n/\ln n$, $D_n/\ln n$ and
$\CL_n/\ln n$.}
\end{center}
\end{figure}

\smallskip\noindent
{\bf Simulations.}
In a first set of computations, two random sequences whose letters are i.i.d.\@
were generated.
On Figure~6,
in the first graph, letters are
equally-likely drawn;
in the second one, they
are drawn with respective probabilities $(p_A,p_C,p_G,p_T)=(0.4,0.3,0.2,0.1)$.
On can visualize the dynamic convergence of $\CL_n/\ln n$, $\ell_n/\ln n$
and of the normalized insertion depth $D_n/\ln n$
(see section~\ref{sec:profondeur}) to their respective constant limits.

Figure~7
is made from simulations of 2,000 random sequences
of length $100,000$ with i.i.d.~letters under the distribution
$(p_A,p_C,p_G,p_T)=(0.6,0.1,0.1,0.2)$.
On the $x$-axis, respectively, lengths of the shortest branches, insertion
depth of the last inserted word, lengths of the longest branches.
On the $y$-axis, number of occurences (histograms).

\begin{proof}[Proof of Lemma~\ref{lemme}]
Since $X_{n}(s)=k$ for $n=T_{k}(s)$ (see Equation~(\ref{duality})), by
monotonicity arguments, it is sufficient to prove that
\[
\frac{\ln T_{k}(s)}{k}\limite{k \to \infty}{a.s.}h(s).
\]
Let  $\varepsilon _r(s)\egaldef Z_r(s)-\E\left[ Z_r(s)\right]$,
so that $T_k(s)$ admits the decomposition
\[
T_k(s)=\E \left[ T_k(s)\right]+\sum _{r=1}^k\varepsilon _r(s).
\]
If $(M_k(s))_k$ is the martingale defined by
\[
M_k(s)\egaldef\sum _{r=1}^k\varepsilon _r(s),
\]
taking the logarithm in the preceding
equation leads to
\begin{equation}
\label{logTk}
\ln T_k(s)=\ln \E\left[ T_k(s)\right]
+\ln\left( 1+\frac{M_k(s)}{\E\left[ T_k(s)\right]}\right) .
\end{equation}

\vskip 5pt\noindent
$\bullet$
It is shown in~\cite{Daudin-Robin} that
$\E\left[ Z_n(s)\right] =1/p\left(s^{(n)}\right)$ so that the sequence
$\frac 1n\ln \E\left[ Z_n(s)\right]$
converges to $h(s)$ as $n$ tends to infinity
($h(s)$ is defined by~(\ref{assumption})).
Since $\E\left[ T_k(s)\right] =\sum _{r=1}^k\E\left[ Z_r(s)\right]$
(see~(\ref{decompcle})), the equality
\[
\lim _{k\to\infty}\frac 1k\ln \E\left[ T_k(s)\right] =h(s)
\]
is a straightforward consequence of the following elementary result:
if $(x_k)_k$ is a sequence of positive numbers such that
$\lim _{k\to\infty}\frac 1k\ln \left( x_k\right) =h>0$, then
$\lim _{k\to\infty}\frac 1k\ln \left( \sum _{r=1}^kx_r\right) =h$.

\vskip 5pt\noindent
$\bullet$
The martingale $(M_k(s))_k$ is square integrable;
its increasing process is denoted by
$\Big(\langle M(s)\rangle _k\Big) _k$.
\cite{Daudin-Robin} have shown that the variance of $Z_r(s)$ satisfies
$\V\left[ Z_r(s)\right]\leq 4r/p\left( s^{(r)}\right) ^2$, so that
\[
\langle M(s)\rangle _k=O\left( ke^{2kh(s)}\right) .
\]
One can thus apply the Law of Large Numbers for martingales
(see~\cite{Duflo} for a reference on the subject):
for any $\alpha >0$,
\[
M_k(s)=O\left(
\langle M(s)\rangle _k^{1/2}
\left( \ln\langle M(s)\rangle _k\right) ^{\frac{1+\alpha}2}
\right)
\hskip 10pt a.s.
\]
Consequently,
\[
\frac{M_k(s)}{\E\left[ T_k(s)\right]}=O\left( k^{1+\alpha /2}\right)
\hskip 10pt a.s.
\]
which completes the proof of Lemma~\ref{lemme}.
\end{proof}

\begin{proof}[Proof of Theorem~\ref{cvps}]
It is inspired from \cite{Pittel}.
Clearly the definition given in Equation (\ref{longueurs}) yields
\[
\ell_{n}\leq X_{n}(s_{+})
\quad \mbox{and} \quad
\CL_{n}\geq X_{n}(s_{-})
\]
(definitions of $s_+$ and $s_-$ were given in~(\ref{h+-})).
Hence, by Lemma~\ref{lemme}
\[
\limsup_{n\to \infty} \frac{\ell_{n}}{\ln n} \leq \frac{1}{h_+},\quad
\liminf_{n\to \infty} \frac{\CL_{n}}{\ln n}\geq \frac{1}{h_-}\quad
\mbox{a.s.}
\]

\noindent
{\it $\bullet$ Proof for $\ell_{n}$}

For any integer $r$,
\begin{equation}
\label{somme}
\PP(\ell_{n}\leq r-1)\leq 
\sum_{s^{(r)}\in{\cal A}^r}\PP(X_{n}(s)\leq r-1)\leq\sum_{s^{(r)}\in{\cal A}^r}
\PP(T_{r}(s)\geq n),
\end{equation}
where the above sums are taken over the set ${\cal A}^r$ of words with
length $r$ (for a proper meaning of this formula, one should replace $s$ by
any infinite word having $s^{(r)}$ as prefix, in both occurences). We abuse of
this notation from now on.
Since the generating functions $\Phi(s^{j},t)$ are defined for any
$1\leq t< \min\{\gamma ^{-1},1+\kappa p(s^{(r)})\}$ and $j \leq r$
(see Assertion i) in Proposition~\ref{generatrice}),
each term of the sum (\ref{somme}) can be controlled by
\[
\PP(T_{r}(s) \geq n)\leq t^{-n}\E[t^{T_{r}(s)}]\leq t^{-n}
\prod_{j=1}^{r}\Phi(s^{(j)},t).
\]
In particular, bounding above all the overlapping functions
$\ind{s_{j}\ldots s_{1}=s_r\dots s_{r-j+1}}$ by $1$ in (\ref{gammadelta}),
we deduce from (\ref{gen}) and from Assertion ii) of
Proposition~\ref{generatrice} that
\[
\PP(T_{r}(s) \geq n)\leq t^{-n}\prod_{j=1}^{r}
\Biggl(1+(1-t)\biggl(
\sum_{\nu=1}^{j}\frac{1}{t^{\nu}p(s^{(\nu)})}+\frac{\kappa'}{1-\gamma t}
\biggr)\Biggr)^{-1}.
\]
Let $0 <\ep <1$. There exists a constant $c_{2}\in ]0,1[$ depending only
on $\ep$ such that
\[
p(s^{(j)})> c_{2} \alpha ^{j},
\quad{\rm with}\quad
\alpha\egaldef\exp (-(1+\ep ^{2})h_+)
\]
(for the sake of brevity $c$, $c_{1}$ and $c_{2}$ denote different constants
all along the text).
We then have
\[
\PP(T_{r}(s)\geq
n) \leq t^{-n}\prod_{j=1}^{r}\Biggl(1+(1-t)\Bigl(\frac{1-(\alpha
  t)^{-j}}{c_{2}(\alpha t -1)}+\frac{\kappa'}{1-\gamma t}\Bigr)\Biggr)^{-1}.
\]
Choosing $t=1+c_{2} \kappa \alpha ^{r}$, Inequality~(\ref{inegGamma}) is valid
if $r$ is large enough, so that
\[
\PP(T_{r}(s) \geq n)\leq ct^{-n}\prod_{j=1}^{r}\Biggl(
1-\kappa \alpha ^{r-j}
\frac{\alpha ^{j}-(1+c_{2} \kappa \alpha ^{r})^{-j}}
{\alpha (1+c_{2} \kappa \alpha ^{r})-1}
-\frac{\alpha ^{r} c_{2}\kappa \kappa'}{1-\gamma (1+c_{2} \kappa \alpha ^{r})}
\Biggr)^{-1}.
\]
Moreover since obvioulsy
\[
\lim_{j \to \infty}
\frac{\alpha ^{j}-(1+c_{2} \kappa \alpha ^{r})^{-j}}
{\alpha (1+c_{2} \kappa \alpha ^{r})-1}
=\frac{1}{1-\alpha},
\]
and
$c_{2}\kappa \kappa'/\bigl(1-\gamma (1+c_{2} \kappa \alpha ^{r})\bigr)$
is uniformly bounded in $r$, there exist two positive constants
$\lambda$ and $L$ independent of $j$ and $r$ such that 
\[
\PP(T_{r}(s) \geq n)\leq
(1+ c_{2} \kappa\alpha ^{r})^{-n}L
\prod_{j=1}^{r}\Bigl(1-\lambda \alpha^{r-j}\Bigr)^{-1}.
\]
In addition, the product can be bounded above by
\[
\prod_{j=1}^{r}\Bigl(1-\lambda \alpha^{r-j}\Bigr)^{-1} \leq
\prod_{j=0}^{\infty}\Bigl(1-\lambda \alpha^{j}\Bigr)^{-1}=R <\infty .
\]
Consequently,
\[
\PP(T_{r}(s) \geq n) \leq L R (1+ c_{2} \kappa \alpha ^{r})^{-n}.
\]
For $r=\lfloor (1-\ep)\frac{\ln n}{h_+}\rfloor$ and $\ep$ small enough, there
exists a constant $R'$ such that
\[
\PP(T_{r}(s) >
n)\leq R' \exp(-c_{2} \kappa n^{\theta}),
\]
where $\theta=\ep-\ep^2+\ep^3 >0$.
We then deduce from (\ref{somme}) that 
\[
\PP(\ell_{n}\leq r-1)\leq 4^{r}R'\exp(-c_{2}\kappa n^{\theta}),
\]
which is the general term of a convergent series.
Borel-Cantelli Lemma applies so that
\[
\liminf_{n\to \infty}
\frac{\ell_{n}}{\ln n}\geq \frac{1}{h_+}\quad \mbox{a.s.}
\]

\noindent{\it $\bullet$ Proof for $\CL_{n}$}

To complete the proof, one needs to show that
\[
\limsup_{n \to \infty} \frac{\CL_{n}}{\ln n}\leq \frac{1}{h_{-}}
\quad\mbox{a.s.}
\]
Again, since $X_{n}(s)=k$ for $n=T_{k}(s)$, by monotonicity arguments
it suffices to show that
\[
\liminf_{k \to \infty} \min_{s^{(k)}\in{{\cal A}^k}}\frac{\ln T_{k}(s)}{k}
\geq h_{-} \quad
\mbox{a.s.}
\]
(notations of~(\ref{somme})).

Let $0<\ep<1$.
As in the previous proof for the shortest branches, it suffices to bound above
\[
\PP\left(\min_{s^{(k)}\in{{\cal A}^k}}T_{k}(s)<e^{kh_{-}(1-\ep)}\right)
\]
by the general term of a convergent series to apply Borel-Cantelli Lemma.
Obviously,
\[
\PP\left( \min_{s^{(k)}\in{{\cal A}^k}}T_{k}(s)<e^{kh_{-}(1-\ep)}\right) \leq
\sum_{s^{(k)}\in{{\cal A}^k}}\PP\left( T_{k}(s)<e^{kh_{-}(1-\ep)}\right).
\]
If $t$ is any real number in $]0,1[$ and if $n\egaldef \exp(kh_{-}(1-\ep))$,
\[
\PP\left( T_{k}(s)<e^{kh_{-}(1-\ep)}\right)
=\PP\left( t^{T_{k}(s)} > t^{n}\right)
\]
and the decomposition (\ref{decompcle}), together with the independence of the
$Z_{r}(s)$ for $1\leq r\leq k$, yield
\[
\PP\left( t^{T_{k}(s)} > t^{n}\right) \leq
t^{-n}\prod_{r=1}^{k}\E\bigl[t^{Z_{r}(s)}\bigr].
\]
The proof consists in bounding above
\begin{equation}\label{finalsum}
\sum_{s^{(k)}\in{{\cal A}^k}}t^{-n}\prod_{r=1}^{k}\E\bigl[t^{Z_{r}(s)}\bigr]
\end{equation}
by the general term of a convergent series, taking $t$ of the form
\[
t\egaldef
(1+c/n)^{-1}
\]
so that the sequence $(t^n)_n$ is bounded.

The generating function of $Z_{r}(s)$ is given by \cite{Daudin-Robin} and
strongly depends on the overlapping structure of the word $s^{(r)}$.
As $0<t<1$, this function is well defined at $t$ and is given by
(see Assertion i) of Proposition~\ref{generatrice})
\begin{equation}
\label{genZj(s)}
\E\bigl[t^{Z_{r}(s)}\bigr]=1-\frac{(1-t)}
{t^{r}p\bigl(s^{(r)}\bigr)\bigl(\gamma_{r}(t)+(1-t)\delta_{r}(t^{-1})\bigr)},
\end{equation}
where $\gamma_{r}(t)$ and $\delta_{r}(t)$ are defined in (\ref{gammadelta}).
Moreover, from Assertion ii) of Proposition~\ref{generatrice}, it is obvious
that there exists a constant $\theta$ independent of $r$ and $s$ such that,
\begin{equation}
\label{majogamma}
\gamma_{r}(t) \leq 1+\theta (1-t).
\end{equation}
Besides, by elementary change of variable, one has successively
\[
\begin{array}{rcl}
t^{r}p\bigl(s^{(r)}\bigr)\delta_{r}(t^{-1})
&=&\DD \sum_{m=1}^{r}\ind{s_r\ldots s_{r-m+1}=s_{m}\ldots s_1}
\frac{t^{r}p\bigl(s^{(r)}\bigr)}{t^{m}p\bigl(s^{(m)}\bigr)}\\
&=&\DD \sum_{m=1}^{r}\ind{s_r\ldots s_{m}=s_{r-m+1}\ldots s_1}t^{m-1}
\frac{p\bigl(s^{(r)}\bigr)}{p\bigl(s^{(r-m+1)}\bigr)}\\
&=&\DD \sum_{m=1}^{r}\ind{s_r\ldots s_{m}=s_{r-m+1}\ldots s_1}t^{m-1}
\frac{p\bigl(s^{(m)}\bigr)}{p\bigl(s_m\bigr)}.
\end{array}
\]
When $m$ is large enough, $h_-$'s definition implies that
\[
p\bigl(s^{(m)}\bigr)\leq\beta^{m},
\quad{\rm where}\quad 
\beta\egaldef\exp(-(1-\ep^{2})h_{-}),
\]
so that there exists positive constants $\rho$ and $c$ such that, for any $r$,
\begin{equation}
\label{truc}
p\bigl(s^{(r)}\bigr)\leq c\beta^{r}{\hskip 10pt\rm and\hskip 12pt}
t^{r}p\bigl(s^{(r)}\bigr)\delta_{r}(t^{-1})\leq
1+\rho\sum_{m=2}^{r}\ind{s_r\ldots s_{m}=s_{r-m+1}\ldots s_1}\beta ^m.
\end{equation}
Thus Formula~(\ref{genZj(s)}) with inequalities~(\ref{majogamma})
and~(\ref{truc}) yield, for any $r\leq k$,
\begin{equation}
\label{majoration}
\E[t^{Z_r(s)}]\leq 1-\frac{1}
{c\beta^r\Bigl(\frac{1}{1-t}+\theta\Bigr)+1+q_{k}(s)},
\end{equation}
where $q_{k}(s)$, that depends on the overlapping structure of $s^{(k)}$, is
defined by
\[
q_{k}(s)\egaldef \rho\max_{2\leq r \leq k} \sum_{m=2}^{r}
\ind{s_r\ldots s_m=s_{r-m+1}\ldots s_1}\beta^{m}.
\]
Note that whatever the overlapping structure is, $q_{k}(s)$ is controlled by
\begin{equation}
\label{qk}
0 \leq q_{k}(s) \leq \frac{\rho}{1-\beta}.
\end{equation}
Thus,
\[
\prod_{r=1}^{k}\E[t^{Z_r(s)}] \leq \exp \biggl[
-\sum_{r=1}^{k} \ln \Bigl(1-\frac{1}
{c\beta^r\bigl((1-t)^{-1}+\theta\bigr)+1+q_{k}(s)}\Bigr)^{-1}
\biggr].
\]
Since the function $x \mapsto \ln 1/(1-x)$ is increasing, comparing
this sum with an integral and after the change of variable
$y=c\beta^x\bigl((1-t)^{-1}+\theta\bigr)$, one obtains
\[
\prod_{r=1}^{k}\E[t^{Z_r(s)}] \leq \exp \biggl[
-\frac{1}{\ln \beta^{-1}}
\int_{c\beta^k\bigl((1-t)^{-1}+\theta\bigr)}^{c((1-t)^{-1}+\theta)}
\ln \Bigl(1-\frac{1}{y+1+q_{k}(s)}\Bigr)^{-1}\frac{dy}{y}
\biggr].
\]
This integral is convergent in a neighbourhood of $+\infty$, hence there exists
a constant $C$, independent of $k$ and $s$ such that 
\begin{equation}
\label{inegproduit}
\prod_{r=1}^{k}\E[t^{Z_r(s)}]\leq C\exp\biggl[
-\frac{1}{\ln \beta^{-1}}
\int_{c\beta^k\bigl((1-t)^{-1}+\theta\bigr)}^{+\infty}
\ln \Bigl(1-\frac{1}{y+1+q_{k}(s)}\Bigr)^{-1}\frac{dy}{y}
\biggr].
\end{equation}
The classical dilogarithm $\dilog_{2}(z)=\sum_{k\geq 1}z^k/k^2$, analytically
continued to the complex plane slit along the ray $[1,+\infty [$, satisfies
$\frac d{dy}\dilog_{2}(-\frac vy)=\frac 1y\log (1+v/y)$.
This leads to the formula
\[
\int_{a_k}^{+\infty}
\ln\left( 1-\frac{1}{y+1+q_{k}(s)}\right) ^{-1}\frac {dy}y
=\dilog_{2}\left( -\frac{q_k(s)}{a_k}\right)
-\dilog_{2}\left(-\frac{1+q_k(s)}{a_k}\right)
\]
with the notation $a_k=c\beta ^k\bigl((1-t)^{-1}+\theta\bigr)$.
Choosing $t=(1+c/n)^{-1}$ yields readily
\begin{equation}
\label{ak}
a_k
\mathop{\sim}_{k\to +\infty}
\exp (-kh_-(\varepsilon-\varepsilon ^2)).
\end{equation}
Moreover, in a neighbourhood of $-\infty$,
\begin{equation}
\label{asymptLi2}
\dilog_{2}(x)=-\frac 12\ln ^2(-x)-\zeta (2)+O(\frac 1{x}),
\end{equation}
and the function $\dilog_{2}(x)+\frac 12\ln ^2(-x)$ is non-decreasing on
$]-\infty ,0[$, so that
\begin{equation}\label{inegLi2}
\left\{
\begin{array}{ll}
\displaystyle
\dilog_{2}(x)\geq -\frac 12\ln ^2(-x)-\zeta (2)\hskip 30pt &(x<0)\\
\displaystyle
\dilog_{2}(x)\leq -\frac 12\ln ^2(-x)-\frac{\zeta (2)}2
\hskip 30pt &(x<-1),
\end{array}
\right.
\end{equation}
noting that $\dilog_{2}(-1)=-\displaystyle\frac{\zeta (2)}2$.
Hence, if $k$ is such that $a_k<1$,
\begin{equation}
\label{majofinal}
\int_{a_k}^{+\infty}
\ln\left( 1-\frac{1}{y+1+q_{k}(s)}\right) ^{-1}\frac {dy}y
\geq
\dilog_{2}\left( -\frac{q_k(s)}{a_k}\right)
+\frac 12\ln ^2(a_k)+\frac{\zeta (2)}2\end{equation}
with $\ln a_k$ being asymptotically proportional to $k$ because of~(\ref{ak}).
Thus, the behaviour of the integral in (\ref{inegproduit}) as $k$ tends to
$+\infty$ depends on the asymptotics of $q_{k}(s)$.

\smallskip
Let $z_k \egaldef \exp\bigl(-\sqrt{k}\bigr)$.
The end of the proof consists, for a given $k$, in splitting the
sum~(\ref{finalsum}) into prefixes $s^{(k)}$ that respectively satisfy
$q_k(s) < \exp\bigl(-\sqrt{k}\bigr)$
or $q_k(s)\geq\exp\bigl(-\sqrt{k}\bigr)$.
These two cases correspond to words that respectly have few or many overlapping
patterns.
The choice $z_k=\exp\bigl(-\sqrt{k}\bigr)$ is arbitrary and many other sequence
could have been taken provided that they converge to zero with a speed
of the form $\exp [-o(k)]$.

\vskip 5pt
First let us consider the case of prefixes $s^{(k)}$ such that 
$q_k(s) < \exp\bigl(-\sqrt{k}\bigr)$.
For such words,~(\ref{inegLi2}) and~(\ref{majofinal}) imply that
\[
\int _{a_k}^{+\infty}\ln
\left( 1-\frac{1}{y+1+q_{k}(s)}\right) ^{-1}\frac {dy}y
\geq
-\frac 12\ln ^2\left(\frac{z_k}{a_k}\right)
+\frac 12\ln ^2(a_k)-\frac{\zeta (2)}2,
\]
the second member of this inequality being, as $k$ tends to infinity,
of the form
$$
k\sqrt kh_-(\ep-\ep ^2)+O(k).
$$
Consequently,
\[
\prod_{r=1}^{k}\E[t^{Z_r(s)}] \leq \exp\left[
  -\frac{\ep}{1+\ep }k^{3/2}+O(k)\right]. 
\]
There are $4^{k}$ words of length $k$, hence very roughly, by taking the sum
over the prefixes $s^{(k)}$ such that $q_{k}(s) < z_{k}$, and since $t^{-n}$
is bounded, the contribution of these prefixes to the sum~(\ref{finalsum})
satisfies
\[
\sum_{s^{(k)}\in{{\cal A}^k},~q_{k}(s) <z_{k}}
t^{-n}\prod_{r=1}^{k}\E\bigl[t^{Z_{r}(s)}\bigr]
\leq 4^{k}\exp\left[ -\frac{\ep}{1+\ep}k^{3/2}+O(k)\right] ,
\]
which is the general term of a convergent series.

\vskip 5pt
It remains to study the case $q_{k}(s) \geq z_{k}$.
For such words, let us only consider the inequalities (\ref{qk})
and~(\ref{inegproduit}) that lead to
\[
\prod_{r=1}^{k}\E[t^{Z_r(s)}] \leq
C\exp \biggl[-\frac{1}{\ln \beta^{-1}}
\int_{a_k}^{+\infty}
\ln \Bigl(1-\frac{1} {y+1+\rho (1-\beta)^{-1}}\Bigr)^{-1}
\frac{dy}{y}\biggr].
\]
Since $x\leq\log(1-x)^{-1}$, after some work of integration, 
\begin{equation}
\label{prodqk>sqrtk}
\prod_{r=1}^{k}\E[t^{Z_r(s)}] \leq
\exp\Bigl(-\frac{\ep}{1+\ep}k+o(k)\Bigr).
\end{equation}
The natural question arising now is: how many words $s^{(k)}$ are there,
such that $q_{k}(s) \geq z_{k}$ ?
Let us define 
\[
E_{k} \egaldef \Big\{s^{(k)},~q_{k}(s) \geq e^{-\sqrt k}\Big\}.
\]
The definition of $q_{k}(s)$ implies clearly that
\[
E_{k}\subseteq \Big\{ s^{(k)},~\exists r \leq k,~\rho
\sum_{m=2}^{r}\ind{s_{r}\ldots s_{m}=s_{r-m+1} \ldots s_{1}}\beta^{m}
\geq e^{-\sqrt k}\Big\} .
\]
For any $r \leq k$ and $x>0$, let us define the set
\[
S_{r}(x)\egaldef \Big\{s^{(k)},~\sum_{m=2}^{r}
\ind{s_{r}\ldots s_{m}=s_{r-m+1} \ldots s_{1}}\beta^{m}<x\Big\}.
\]
For any $l\in\{ 2,\dots ,r\}$, one has the following inclusion
\[
\bigcap_{ m=2}^{\ell}\Big\{s^{(k)},~
\ind{s_{r}\ldots s_{m}=s_{r-m+1} \ldots s_{1}}=0\Big\}
\subset
S_{r}\Bigl(\frac{\beta^{\ell +1}}{1-\beta}\Bigr).
\]
If the notation $B^{c}$ denotes the complementary set of $B$ in $\CA^{k}$,
\[
S_{r}^c\Bigl(\frac{\beta^{\ell +1}}{1-\beta}\Bigr)
\subset
\bigcup_{ m=2}^{\ell}\Big\{ s^{(k)},~
\ind{s_{r}\ldots s_{m}=s_{r-m+1} \ldots s_{1}}=1\Big\}
\]
Since $e^{-\sqrt k}=\rho\beta^{\ell +1}(1-\beta )^{-1}$ for
$\ell \egaldef\sqrt{k}/\ln\left(\beta^{-1}\right)
+\ln\left( \rho ^{-1}(1-\beta )\right) /\ln\beta$,
\[
E_{k} \subset \bigcup_{r=1}^{k}\bigcup_{m=2}^{\lfloor\ell\rfloor +1}
\Big\{s^{(k)},~
\ind{s_{r}\ldots s_{m}=s_{r-m+1} \ldots s_{1}}=1\Big\},
\]
so that the number of words $s^{(k)}$ such that $q_{k}(s) \geq z_{k}$ is
bounded above by
\begin{equation}
\label{cardEk}
\# E_{k}\leq \sum_{r=1}^{k}\sum_{m=2}^{\lfloor\ell\rfloor +1}4^{m-1}
\in O(k4^{\sqrt k /\ln\bigl(\beta^{-1}\bigr)}.
\end{equation}
Putting~(\ref{prodqk>sqrtk}) and ~(\ref{cardEk}) together is sufficient to
show that the contribution of prefixes $s^{(k)}$ such that $q_k(s)\geq z_k$
to the sum~(\ref{finalsum}), namely
\[
\sum_{s^{(k)}\in{{\cal A}^k},~q_{k}(s)\geq z_{k}}
t^{-n}\prod_{r=1}^{k}\E\bigl[t^{Z_{r}(s)}\bigr],
\]
is the general term of a convergent series too.

Finally, the whole sum~(\ref{finalsum}) is the general term of a convergent
series, which completes the proof of the inequality
\[
\limsup_{n \to \infty}\frac{\CL_{n}}{\ln n}\leq \frac{1}{h_{-}}
\quad \mbox{a.s.}\]
\end{proof}
\section{Insertion depth}\label{sec:profondeur}

This section is devoted to the asymptotic behaviour of the insertion depth
denoted by $D_{n}$ and to the length of a path randomly and uniformly chosen
denoted by $M_{n}$ (see section~\ref{sec:def}).
$D_{n}$ is defined as the length of the path leading to
the node where $W(n)$ is inserted.
In other words, $D_{n}$ is the amount of digits to be checked before the
position of $W(n)$ is found.
Theorem~\ref{cvps} immediately implies a first asymptotic result on
$D_{n}$. Indeed, $D_{n}=\ell_{n}$ whenever $\ell_{n+1}>\ell_{n}$, which
happens infinitely often a.s., since  $\lim_{n \to \infty} \ell_{n}=\infty$
a.s.
Hence, 
\[
\liminf_{n \to \infty}\frac{D_{n}}{\ln n}
=\liminf_{n \to\infty}\frac{\ell_{n}}{\ln n}=\frac{1}{h_+}\quad \mbox{a.s.}
\]
Similarly, $D_{n}=\CL_{n}$ whenever $\CL_{n+1}>\CL_{n}$, and hence 
\[
\limsup_{n \to \infty}\frac{D_{n}}{\ln n}
=\limsup_{n \to \infty}\frac{\CL_{n}}{\ln n}=\frac{1}{h_-}\quad \mbox{a.s.}
\]
Theorem~\ref{thm} states full convergence in probability of these random
variables to the constant $1/h$.
\begin{thrm}
\label{thm}
\[\frac{D_{n}}{\ln n} \limite{n \to \infty}{P}\frac{1}{h}\quad\mbox{and} \quad
\lim_{n\to \infty}\frac{M_{n}}{\ln
  n} \limite{n \to \infty}{P}\frac{1}{h}.\]
\end{thrm}
\begin{rmrk}
For an i.i.d.\@ sequence $U=U_1 U_2 \ldots$, in the case when the random
variables $U_{i}$ are not uniformly distributed in $\{A, C, G, T\}$,
Theorem~\ref{thm} implies that $\frac{D_{n}}{\ln n}$ does not converge
a.s. because
\[
\limsup_{n\to \infty}\frac{D_{n}}{\ln n}
\geq \frac{1}{h}>\frac{1}{h_{+}}=\liminf_{n\to \infty}\frac{D_{n}}{\ln n}.
\]
\end{rmrk}
\begin{proof}[Proof of Theorem~\ref{thm}]
It suffices to consider $D_{n}$ since, by definition of $M_{n}$,
\[
\PP(M_{n}=r)=\frac{1}{n}\sum_{\nu=1}^{n}\PP(D_{\nu}=r).
\]
Let $\ep >0$.
To prove Theorem~\ref{thm}, we get the convergence
$\lim_{n \to \infty}\PP(A_{n})=0$, where
\[
A_{n}\egaldef \Big\{ U \in \CA^{\xN},~
\Big|\frac{D_{n}}{\ln n}-\frac{1}{h}\Big|\geq \frac{\ep}{h}\Big\},
\]
by using the obvious decomposition
\[
\PP (A_{n})=\PP\left( \frac{D_{n}}{\ln n}\geq \frac{1+\ep}{h}\right)
+\PP\left( \frac{D_{n}}{\ln n}\leq \frac{1-\ep}{h}\right) .
\]

\noindent
$\bullet$
Because of $X_n$'s definition~(\ref{defXn}),
\[
D_{n}=X_{n-1}\bigl(W(n)\bigr) +1
\]
so that the duality (\ref{duality}) between $X_n(s)$ and $T_k(s)$ implies that
\begin{equation}
\label{decompDn}
\PP\left( \frac{D_{n}}{\ln n}\geq \frac{1+\ep}{h}\right)\leq
\PP\Big( X_{n-1}\bigl(W(n)\bigr)\geq k-1\Big)\leq
\PP\Big( T_{k-1}\bigl(W(n)\bigr)\leq n-1\Big)
\end{equation}
with $k\egaldef\lfloor\frac{1+\ep}{h}\ln n\rfloor$.
Furthermore, 
\[
\PP\Big( T_{k-1}\bigl(W(n)\bigr)\leq n-1\Big)\leq
\PP\Big( \left\{ T_{k-1}\bigl(W(n)\bigr)\leq n-1\right\}\cap B_{n,k_0}\Big)
+\PP\Big( B_{n,k_0}^c\Big)
\]
where $B_{n,k_0}$ is defined, for any $k_0\leq n$, by
\[
B_{n,k_0}\egaldef \bigcap_{k_0\leq j\leq n}\Big\{ U\in \CA^{\xN},~
\Big|\frac{1}{j}\ln \Bigl(\frac{1}{p\bigl(W(n)^{(j)}\bigr)}\Bigr)-h\Big|
\leq\ep^{2}h\Big\}.
\]
Since the sequence $U$ is stationary,
$\PP\big( W(n)^{(j)}\big) =\PP\big( U^{(j)}\big)$
so that Ergodic Theorem implies
\[
\lim _{j\to\infty }
\frac 1j\ln\left( \frac 1{p\big( W(n)^{(j)}\big) }\right)
=h
\hskip 20pt{\rm a.s.}
\]
which leads to $\PP\big( B_{n,k_0}\big) =1$ when both $k_0$ and $n$
are large enough.
If $\CS_{n,k_0}$ denotes the set of words
\[
\CS _{n,k_0}\egaldef
\left\{
s^{(n)}\in\CA ^n,~\forall j\in\{ k_0,\dots ,n\}~
\Big|\frac 1j\ln\left(\frac 1{p\bigl( s^{(j)}\bigr)}\right) -h\Big|\leq\ep ^2h
\right\},
\]
when $k_0$ and $n$ are large enough,
\[
\begin{array}{rcl}
\PP\Big( T_{k-1}\bigl(W(n)\bigr)\leq n-1\Big)&\leq &
\DD\sum _{s^{(n)}\in\CS _{n,k_0}}
\PP\Big( W(n)^{(n)}=s^{(n)},~T_{k-1}(s)\leq n-1\Big)\\
&\leq &
\DD\sum _{s^{(n)}\in\CS _{n,k_0}}\PP\Big( T_{k-1}(s)\leq n-1\Big) .
\end{array}
\]
Such a probability has already been bounded above at the end of
Theorem~\ref{cvps}'s proof;
similarly,
\begin{equation}
\label{majoSkk0}
\sum _{s^{(n)}\in\CS _{n,k_0}}\PP\Big( T_{k-1}(s)\leq n-1\Big)
=O\left(
n\exp\Big(-\frac{\ep}{1+\ep}n+\frac{\ln 4}{(1-\ep ^2)h}\sqrt n\Big)
\right)
\end{equation}
so that~(\ref{decompDn}) and~(\ref{majoSkk0}) show that
$\PP\left( \frac{D_{n}}{\ln n}\geq \frac{1+\ep}{h}\right)$
tends to zero when $n$ goes off to infinity.

\vskip 5pt\noindent
$\bullet$
Our argument showing that
$\PP\left( \frac{D_{n}}{\ln n}\leq \frac{1-\ep}{h}\right)$
tends to zero when $n$ tends to infinity is similar.
If now
$k\egaldef\lfloor\frac{1-\ep}{h}\ln n\rfloor$,
\[
\PP\left( \frac{D_{n}}{\ln n}\leq \frac{1-\ep}{h}\right)\leq
\PP\Big( X_{n-1}\bigl(W(n)\bigr)\leq k-1\Big) =
\PP\Big( T_{k}\bigl(W(n)\bigr)\geq n\Big) ,
\]
so that
\[
\PP\left( \frac{D_{n}}{\ln n}\leq \frac{1-\ep}{h}\right)\leq
\PP\Big( \left\{ T_{k}\bigl( W(n)\bigr)\geq n\right\}\cap B_{n,k_0}\Big)
+\PP\big( B_{n,k_0}^c\big).
\]
As before, $\PP\big( B_{n,k_0}^c\big)=0$ when $k_0$ and $n$ are large enough
and
\[
\begin{array}{rcl}
\PP\Big( T_{k}\bigl(W(n)\bigr)\geq n\Big)&\leq &
\DD\sum _{s^{(n)}\in\CS _{n,k_0}}
\PP\Big( W(n)^{(n)}=s^{(n)},~T_{k}(s)\geq n\Big)\\
&\leq &
\DD\sum _{s^{(n)}\in\CS _{n,k_0}}\PP\Big( T_{k}(s)\geq n\Big) .
\end{array}
\]
Like in the proof of Theorem~\ref{cvps}, on shows that
\[
\sum _{s^{(n)}\in\CS _{n,k_0}}\PP\Big( T_{k}(s)\geq n\Big)
=O\left(
4^n\exp\left( -\kappa n^\theta /2\right)\right)
\]
which implies that 
$\PP\left( \frac{D_{n}}{\ln n}\leq \frac{1-\ep}{h}\right)$
tends to zero when $n$ tends to infinity.
The proof of Theorem~\ref{thm} is complete.
\end{proof}
\newpage
\appendix
\section{Domain of definition of the generating function $\Phi(s^{(r)},t)$}
\label{sec:A}
\subsection{Proof of Assertion ii)}
There exists a function $K(s_{1},s_{r},m)$ uniformly bounded by the constant
$$
K\egaldef \sup_{s_{1},s_{r},m}|K(s_{1},s_{r},m))|
$$
such that
\begin{equation}
\label{stationnaire}
Q ^{m}(s_{1}, s_{r})-p\bigl(s_{r}\bigr)=K(s_{1},s_{r},m)\gamma ^{m}
,\end{equation}
where $\gamma$ is the second eigenvalue of the transition
matrix. Consequently,
\begin{eqnarray*}
|\gamma_{r}(t)-1| &=&\Bigl|\frac{1-t}{tp\bigl(s_{r}\bigr)}\sum_{m \geq
  1}K(s_{1},s_{r},m)(\gamma t)^{m}\Bigr| \\
&\leq& \frac{\gamma K}{\min_{u}p(u)}\frac{|1-t|}{1-\gamma |t|}.
\end{eqnarray*}
Hence Assertion ii) holds with $\kappa '\egaldef \gamma K/\min_{u}p(u)$.


\subsection{Proof of Assertion i)}
On the unit disc $|t|<1$, the series
\begin{equation}
\label{serie}
S(t)\egaldef\frac{1}{t}\sum_{m\geq 1}Q ^{m}(s_{1}, s_{r})t^{m}
\end{equation}
is convergent and one has the decomposition
\[\frac{1-t}{p\bigl(s_{r}\bigr)t} \sum_{m\geq 1}Q ^{m}(s_{1}, s_{r})t^{m}=1
+\frac{1-t}{p\bigl(s_{r}\bigr)t} \sum_{m\geq 1}\bigl[Q ^{m}(s_{1},
s_{r})-p\bigl(s_{r}\bigr)\bigr]t^{m}.\]
The function
\[
\sum_{m\geq 1}\bigl[Q ^{m}(s_{1}, s_{r})-p\bigl(s_{r}\bigr)\bigr]t^{m}
\]
is analytically continuable to the domain $\gamma |t|<1$, and then the series
\[
\frac{1-t}{t p\bigl(s_{r}\bigr)}\sum_{m\geq 1}Q ^{m}(s_{1}, s_{r})t^{m}
\]
converges on the same domain.
One has to determine the zeroes of
\begin{eqnarray*}
D(t)\egaldef p\bigl(s^{(r)}\bigr)t^{r}&
+&\frac{(1-t)p\bigl(s^{(r)}\bigr)t^{r}}{p\bigl(s_{r}\bigr)t}
\sum_{z\geq 1}t^{z}\bigl[Q^{z}(s_{1},s_{r})-p\bigl(s_{r}\bigr)\bigr]\\
&+&(1-t)\bigl[1+ \sum_{j=2}^{r}t^{j-1}\frac{p\bigl(s^{(j)}\bigr)}{p(s_j)}
\ind{s_{r}\ldots s_{j}=s_{r-j+1}\ldots s_{1}}\bigr].
\end{eqnarray*}
Assuming that some $0<t<1$ were a real root of $D(t)$, then
\begin{eqnarray*}
0&<&\frac{(1-t)p\bigl(s^{(r)}\bigr)t^{r}}{p\bigl(s_{r}\bigr)t}\sum_{z\geq
  1}t^{z}Q^{z}(s_{1},s_{r})\\
&=&(t-1)\bigl[1+ \sum_{j=2}^{r}t^{j-1}\frac{p\bigl(s^{(j)}\bigr)}{p(s_j)}
\ind{s_{r}\ldots s_{j}=s_{r-j+1}\ldots s_{1}}\bigr] <0.
\end{eqnarray*}
It is thus obvious that there are no real root of $D(t)$ in $]0,1[$.
Moreover, one can readily check that $0$ and  $1$ are not zeroes of $D(t)$.
We now look for a root of the form $t=1+\ep$ with $\ep >0$.
Such an $\ep$ satisfies
\[
\ep =\frac
{\DD (1+\ep)^{r}p\bigl(s^{(r)}\bigr)\Bigl(1-\frac{\ep}{p(s_{r})(1+\ep)}
\sum_{z\geq 1} t^{z}[Q^{z}(s_{1},s_{r})-p(s_{r})]\Bigr)}
{\DD 1+ \sum_{j=2}^{r}(1+\ep )^{j-1}\frac{p\bigl(s^{(j)}\bigr)}{p(s_j)}
\ind{s_{r}\ldots s_{j}=s_{r-j+1}\ldots s_{1}}},
\]
so that
\[
\ep\geq
\kappa p\bigl(s^{(r)}\bigr).
\]
This implies that $\Phi(s^{(r)},t)$ is at least defined on
$\bigl[0, 1+\kappa p\bigl(s^{(r)}\bigr) \bigr[$.
This implies the result.

\bibliographystyle{plainnat}
\bibliography{Quad}
\end{document}